\documentclass{amsart}
\usepackage{amsxtra,amssymb,multicol,graphicx}

\newtheorem{thm}{Theorem}

\newtheorem{question}[thm]{Question}

\newtheorem{cor}[thm]{Corollary}

\newtheorem{conj}[thm]{Conjecture}

\theoremstyle{definition}
\newtheorem{defn}[thm]{Definition}
\newtheorem{tempdef}[thm]{Temporary Definition}

\newtheorem{say}[thm]{}
\newtheorem{exmp}[thm]{Example}

\newtheorem{exrc}[thm]{Exercise}
\newtheorem{prob}[thm]{Problem}
\newtheorem{ques}[thm]{Question}    

\newtheorem{rem}[thm]{Remark}          

\newtheorem{ack}{Acknowledgments}        
   
\newtheorem{warning}[thm]{Warning}  
\newtheorem{defn-thm}[thm]{Definition--Theorem}  
\newtheorem{defn-lem}[thm]{Definition--Lemma}  

\theoremstyle{remark}

\renewcommand{\c}[0]{{\mathbb C}}  

\renewcommand{\o}[0]{{\mathcal O}} 
\newcommand{\z}[0]{{\mathbb Z}}

\renewcommand{\r}[0]{{\mathbb R}} 

\renewcommand{\a}[0]{{\mathbb A}}

\newcommand{\p}[0]{{\mathbb P}}

\newcommand{\q}[0]{{\mathbb Q}}
\newcommand{\map}[0]{\dasharrow}
\newcommand{\qtq}[1]{\quad\mbox{#1}\quad}

\newcommand{\pic}[0]{\operatorname{Pic}}

\newcommand{\cont}[0]{\operatorname{cont}}

\newcommand{\mult}[0]{\operatorname{mult}}

\newcommand{\supp}[0]{\operatorname{Supp}}    
\newcommand{\red}[0]{\operatorname{red}}    
\newcommand{\codim}[0]{\operatorname{codim}}    
    
\newcommand{\proj}[0]{\operatorname{Proj}}

\newcommand{\cent}[0]{\operatorname{center}}
\newcommand{\bs}[0]{\operatorname{Bs}}

\newcommand{\sing}[0]{\operatorname{Sing}}    
\newcommand{\ex}[0]{\operatorname{Ex}}

\newcommand{\cl}[0]{\operatorname{Cl}}

\newcommand{\nec}[0]{\overline{\operatorname{NE}}}

\newcommand{\rup}[1]{\lceil{#1}\rceil}
\newcommand{\rdown}[1]{\lfloor{#1}\rfloor}

\newcommand{\simq}[0]{\sim_{\q}}
\newcommand{\simr}[0]{\sim_{\r}}

\newcommand{\cdiv}[0]{\operatorname{CDiv}}
\newcommand{\wdiv}[0]{\operatorname{Div}}

\newcommand{\tsum}[0]{\operatorname{\textstyle{\sum}}}

\begin{document}
\bibliographystyle{amsalpha}

\title{Exercises in the birational geometry of algebraic varieties}
\author{J\'anos Koll\'ar}

\maketitle


The book \cite{km-book} gave an introduction to
the birational geometry of algebraic varieties,
as the subject stood in 1998.  The developments of the
last decade made the more advanced parts of Chapters
6 and 7 less important and the detailed treatment of
surface singularities in Chapter 4 less necessary.
However, the main parts, Chapters 1--3 and 5, 
 still form the  foundations of the subject.

These notes provide additional exercises to
\cite{km-book}. The main definitions and theorems are recalled
but not proved here. The emphasis is on the many examples
that illustrate the methods, their shortcomings and
some applications.

\section{Birational classification of algebraic surfaces}

For more detail, see \cite{bpv}.
\medskip

The theory of algebraic surfaces rests on the following
three theorems.

\begin{thm} \label{s1.thm}
Any birational morphism between smooth projective surfaces
is a composite of  blow-downs to points.
Any birational map between smooth projective surfaces
is a composite of blow-ups and blow-downs.
\end{thm}

\begin{thm}\label{s2.thm}
 There are 3 species of ``pure-bred''   surfaces:
\begin{description}
\item[\rm (Rational)] For these surfaces the
internal birational geometry is very complicated, but,
up to birational equivalence, we have only $\p^2$. 
These frequently appear in the classical literature
and in ``true'' applications.
\item[\rm (Calabi-Yau)] These are completely classified
(Abelian, K3, Enriques, hyperelliptic) and their geometry
is rich. They are of great interest
to other mathematicians.
\item[\rm (General type)] They have a canonical model
with Du Val singularities and ample canonical class.
Although singular, this is the ``best'' model to work with.
There are lots of these but they appear less frequently
outside algebraic geometry.
\end{description}
There are also two types of ``mongrels'':
\begin{description}
\item[\rm (Ruled)] Birational to $\p^1\times (\mbox{curve of genus $\geq 1$})$.
\item[\rm (Elliptic)]  These fiber over a curve with
general fiber an elliptic curve.
\end{description}
\end{thm}

The ``mongrels'' are usually studied as an afterthought,
 with suitable modifications
of the existing methods. In a general survey, it is best to ignore them.

\begin{thm}\label{s3.thm}
 Assume that $S$ is neither rational nor ruled.
Then there is a unique smooth projective surface $S^{min}$
 birational to $S$
such that every birational map
$S'\map S^{min}$ from a 
smooth projective surface  $S'$ is automatically
a morphism.
\end{thm}

Some of these theorems are relatively easy, and some
condense a long and hard story into a short statement.

\begin{verse}
{\it  The first aim of higher dimensional algebraic geometry is to}\\
{\it  generalize these theorems to dimensions three and up.}
\end{verse}

In these notes we focus only on certain aspects of this
project. Let us start with mentioning the parts that we
will not cover.

The correct higher dimensional analogs of rational surfaces are
{\it rationally connected varieties} and the ruled surfaces
are replaced by {\it rationally connected fibrations.}
We do not deal with them here.
See \cite{k-bull} for an introduction and
\cite{rc-book} for a detailed treatment.

The study of higher dimensional Calabi-Yau varieties is very active,
with most of the effort going into understanding mirror symmetry
rather than developing a general classification scheme.

It is  known that 
any birational {\em map} between smooth projective varieties
is a composite of blow-ups and blow-downs of smooth subvarieties
\cite{wlod, akmw}.  While it is very useful to stay with
these easy-to-understand elementary steps, in practice it 
is very hard to keep track of geometric properties during
  blow-ups. It is much more useful to factor every 
birational {\em morphism} between smooth projective varieties
as a composite of elementary steps. It turns out that
smooth blow-ups do not work (\ref{nonfactor.exrc}).
  From our current point of view,
the natural question is to work with varieties with
terminal singularities and consider the factorization 
of birational morphisms as a special
case of the MMP. However, the following intriguing
problem is still open.

\begin{ques} Let $f:X\to Y$ be a birational morphism between
smooth projective 3-folds. Is $f$ a composite of smooth 
blow-downs and flops?
\end{ques}

\section{Naive minimal models}

This is a more technical version of my notes
\cite{what-mm}.
\medskip

Much of the power of affine algebraic geometry
rests on the basic correspondence
$$
\begin{array}{ccc}
&   \mbox{ring of regular functions} & \\
\Bigl\{
\begin{array}{c}
\mbox{affine}\\
\mbox{schemes}
\end{array}
\Bigr\}
& \stackrel{\textstyle{\longrightarrow}}{\longleftarrow} &
\Bigl\{
\begin{array}{c}
\mbox{commutative}\\
\mbox{rings}
\end{array}
\Bigr\}
\\[-1ex]
& \mbox{spectrum} &
\end{array}
$$
Thus every affine variety is the natural existence domain
for the ring of all regular functions on it.

\begin{exrc} \label{affine.exrc}
Let $X$ be a $\c$-variety of finite type. Prove that $X$ is affine
iff the following 2 conditions are satisfied:
\begin{enumerate}
\item (Point separation) 
For any
 two  points $p \neq q\in X$ there is a regular
function $f$ on $X$ such that $f(p)\neq f(q)$.
\item (Maximality of domain) For any sequence of points
$p_i\in X$ that does not converge to a limit in $X$,
 there is a regular
function $f$ on $X$ such that $\lim f(p_i)$ does not exist.
\end{enumerate}
\end{exrc}

\begin{exrc} Reformulate and prove Exercise \ref{affine.exrc}
for varieties over arbitrary fields.
\end{exrc}

\begin{say}
As we move to more general varieties, this nice correspondence breaks down
in two distinct ways.

{\it Quasi affine varieties.} Let $X:=\a^n\setminus(\mbox{point})$
for some $n\geq 2$. Check that every regular function on
$X$ extends to a regular function on $\a^n$. Thus the 
function theory of $X$ is  rich but the 
natural existence domain
for the ring of all regular functions on $X$ is the larger space
$\a^n$. 
Similarly, if 
$$
X=(\mbox{irreducible affine variety})\setminus
(\mbox{codimension $\geq 2$ subvariety}),
$$
then every regular function on
$X$ extends to a regular function on the
irreducible affine variety.

{\it Projective varieties.} On a projective variety every
regular (or holomorphic)
 function is constant, hence the regular (or holomorphic)
function theory of a projective variety is not interesting.

 On the other hand, a projective variety has 
 many interesting  {\it rational functions}.
That is, functions that can locally be written as the quotient
of two regular
functions. At a point the value of a rational function $f$
can be finite, infinite or undefined. 
 The set of points where
$f$ is undefined  has   codimension $\geq 2$.
 This makes it  hard to
control what happens in codimensions $\geq 2$.

Rational functions on a $k$-variety $X$ form a field  $k(X)$, called the
{\it function field} of $X$.
\end{say}

\begin{exrc} Let $X=(xy-uv=0)\subset \a^4$ and $f=x/u$.
Show that $X$ is normal and $f$ is undefined only at the origin $(0,0,0,0)$.
\end{exrc}

\begin{exrc} Let $X$ be a normal, proper variety
over an algebraically closed field $k$. Prove that
$X$ is projective iff  for any two
 points $p\neq q\in X$ and finite subset $R\subset X$, there is a rational
function $f$ on $X$ such that $f(p)\neq f(q)$ and 
 $f$ is defined at all points of $R$.
\end{exrc}

Following the example of affine varieties we ask:

\begin{question} How tight is the connection between $X$ and $k(X)$?
\end{question}

 Assume that we have
$X_1\subset \p^r$ with coordinates $(x_0:\dots : x_r)$,
$X_2\subset \p^s$ with coordinates $(y_0:\dots : y_s) $
and an isomorphism $\Psi:k(X_1)\cong k(X_2)$.
Then $\phi_i:=\Psi(x_i/x_0)$   are rational functions 
on $X_2$ 
and $\phi^{(-1)}_j:=\Psi^{-1}(y_j/y_0)$  are rational functions 
on $X_1$. 
(Note that $\phi^{(-1)}_j$ is not the inverse of  $\phi_j$.)
Moreover, 
$$
\Phi:(y_0:\dots : y_s) \mapsto 
\bigl(1:\phi_{1}(y_0:\dots : y_s):\cdots : \phi_{r}(y_0:\dots : y_s)\bigr)
$$
defines a rational map $\Phi:X_2\map X_1$ and
$$
\Phi^{-1}:(x_0:\dots : x_r) \mapsto 
\bigl(1: \phi^{(-1)}_{1}(x_0:\dots : x_r):\cdots : 
\phi^{(-1)}_{s}(x_0:\dots : x_r)\bigr)
$$
defines a rational map $\Phi^{-1}:X_2\map X_1$
such that 
$\Psi$ is  induced by pulling back functions by
$\Phi$ and  $\Psi^{-1}$ is  induced by pulling back functions by
$\Phi^{-1}$. 
That is, $X_1$ and $X_2$ are
{\it birational} to each other.

\begin{exrc} Let $C_1, C_2$ be 1-dimensional, irreducible, projective
with all local rings regular. Prove that every birational map
$C_1\map C_2$ is an isomorphism.
\end{exrc}

The situation is more complicated in higher dimensions.
A map with an  inverse is usually an isomorphism, but
this fails in the birational case
 since  $\Phi$ and $\Phi^{-1}$ are not everywhere defined.
The simplest examples are blow-ups and blow-downs.

\begin{say}[Blow-ups]
Let $X$ be a smooth projective variety and $Z\subset X$ a smooth subvariety.
Let $B_ZX$ denote the {\it blow-up} of $X$ along $Z$
and $E_Z\subset B_ZX$ the exceptional divisor.  
We refer to $\pi:B_ZX\to X$ as a blow-up if we imagine that
$B_ZX$ is created from $X$, and a {\it blow-down} if we
start with $B_ZX$  and construct $X$ later.
Note that
 $E_Z$ has codimension 1 and $Z$ has codimension $\geq 2$.
Thus a blow-down decreases the Picard number by 1.

By blowing up repeatedly, starting with any $X$ we can create
more and more complicated varieties with the same function field.
Thus, for a given function field
$K=k(X)$, there is no  ``maximal domain''
where all elements of $K$ are rational functions.
(The inverse limit of all varieties birational to $X$
appears in the literature occasionally as  such a ``maximal domain,'' 
but so far with limited success.)
 On the other hand,  one can look
for a ``minimal domain'' or ``minimal model.''

As a first approximation, a variety $X$ is a
minimal  model if the  underlying space $X$
is the ``best match'' to
the rational
function theory of $X$. 
\end{say}

\begin{exmp} Let $S$ be a smooth projective surface
which is neither rational nor ruled. 
Explain why it makes sense to say that  
 $S^{min}$  (as in  (\ref{s3.thm})) is a  ``minimal domain'' for 
the field $k(S)$.
\end{exmp}

\begin{exrc} Let $X$ be a  projective variety 
that admits a 
 finite morphism to an Abelian variety.
Prove that every rational map $f:Y\map X$
from a smooth projective variety $Y$ to $X$ is
 a morphism.

Thus, if $X$ is smooth, 
 it makes sense to say that $X$ is a  ``minimal domain'' of 
its function field $k(X)$.
\end{exrc}

Not all varieties have a ``minimal domain''
with the above strong properties.

\begin{exmp}
Let ${\mathbf Q}^3\subset \c\p^4$ be the quadric hypersurface 
 given by the equation $x^2+y^2+z^2+t^2=u^2$.
Let 
$$
\pi: (x:y:z:t:u)\map (x:y:z:u-t)
$$
 be the projection 
from the north pole $(0:0:0:1:1)$ to the equatorial plane $(t=0)$. 
Its inverse $\pi^{-1}$ is given by 
$$
(x:y:z:u)\map (2xu: 2yu:2zu:x^2+y^2+z^2-u^2: x^2+y^2+z^2+u^2).
$$
These maps show that the meromorphic function theory of
${\mathbf Q}^3$ is the same as that of $\c\p^3$.

Show that  $\pi$ contracts the lines
$(a\lambda:b\lambda:c\lambda:1:1)$ to the points $(a:b:c:0)$
whenever $a^2+b^2+c^2=0$, 
and $\pi^{-1}$ contracts
the plane at infinity $(u=0)$ to the point $(0:0:0:1:1)$.
Write $\pi$ as a composite of
blow ups and blow downs with smooth centers.

On the other hand, ${\mathbf Q}^3$ and  $\c\p^3$ 
are quite different as manifolds. Show that they have the
same Betti numbers but they are not homeomorphic.
Prove that  ${\mathbf Q}^3$ and $\c\p^3$ both have
Picard number 1.
\end{exmp}

A more subtle example is the following.

\begin{exrc}\label{first.flop.exrc} Let
$Y$ be a smooth projective variety of dimension $3$
and $f,g,h$  general sections of a very ample line bundle $L$ on $Y$.
 Consider the hypersurface 
$$
X:=(s^2f+2stg+t^2h=0)\subset Y\times \p^1_{s:t}.
$$
Show that $X$ is smooth and compute its canonical class.

Show that the projection $\pi:X\to Y$ has degree 2;
let $\tau:X\map X$ be the corresponding Galois involution.
Write it down explicitly in coordinates and decide
 where $\tau$ is  regular.

Show that $X$ contains $(L^3)$ curves of the form
$(\mbox{point})\times \p^1$ and they are
numerically equivalent to each other.
(This may need the Lefschetz theorem 
on the Picard groups of  hyperplane sections.)

Assume that  $Y$
 admits a 
 finite morphism to an Abelian variety. Prove that the folloing hold:
\begin{enumerate}
\item 
 Any smooth projective variety $X'$ that is birational to
$X$ has Picard number at least $\rho(X)$.
\item  If
$X'$ has Picard number $\rho(X)$ then it is isomorphic to $X$.
\item If $(L^3)>1$ then there are nonprojective 
compact complex manifolds $Z$ that are bimeromorphic to
$X$, have Picard number  $\rho(X)$, but are not isomorphic to $X$.
\end{enumerate}
\end{exrc}

\begin{exrc}  Let $X$ be a smooth projective variety
such that $K_X$ is nef. 
Let $f:X\map X'$ be a birational map to a smooth projective variety.
Prove that the exceptional set $\ex(f)$ has codimension
$\geq 2$ in $X$. 
Generalize to the case when $X$ is canonical and $X'$ is terminal
(\ref{7.9}).

Hint: You should find (\ref{exc.eff.exrc}) helpful.
\end{exrc}

\begin{defn} We say that a birational map $f:X_1\map X_2$
{\it contracts} a divisor $D\subset X_1$ if
$f$ is defined at the generic point of $D$ and
$f(D)\subset X_2$ has codimension $\geq 2$.
The map  $f$ is called a  {\it  birational contraction} if  
$f^{-1}$ does not contract any divisor.

A birational map $f:X_1\map X_2$ is called {\it small}
if neither $f$ nor $f^{-1}$ contracts any divisor.

The simplest examples of birational contractions are composites of blow-downs,
but there are many, more complicated, examples.
\end{defn}

\begin{exrc}  Let $f:S_1\map S_2$ be a birational contraction between
smooth projective surfaces. Show that $f$ is a morphism.
\end{exrc}

\begin{exrc}  Let $L,M\subset \p^3$ two  lines intersecting at a point.
 The identity on $\p^3$
induces a rational map
$g:B_LB_M\p^3\map B_L\p^3$. 
(With a slight abuse of notation, we also denote by $L$
the birational transform of $L$ on $B_M\p^3$, etc.)
Show that $g$ is a contraction but  it is
not a morphism. Describe how to factor $g$ into a composite of
smooth blow ups and blow downs.\end{exrc}


There is essentially only one way to write
a birational morphism between smooth surfaces as  a composite
of point blow ups. The next exercise shows that this
no longer holds for 3-folds. 

\begin{exrc}  Let $p\in L\subset \p^3$ be a point on a line.
 Let $C\subset B_L\p^3$ be the preimage
of $p$. Show that the identity on $\p^3$
induces an isomorphism
$B_CB_L\p^3\cong B_LB_p\p^3$. 
\end{exrc}

The next exercise shows that not every 
birational morphism between smooth 3-folds is  a composite
of smooth blow-ups.

\begin{exrc} \label{nonfactor.exrc}
Let $C\subset \p^3$ be an irreducible curve with a unique
singular point which is either a node or a cusp.
Show that $B_C\p^3$ has a unique singular point; call it $p$.
Check that $B_pB_C\p^3$ is smooth.
Prove that  $\pi: B_pB_C\p^3\to \p^3$ can not be written as a composite
of smooth blow-ups.

Write $\pi$ as a composite
of two smooth blow-ups and a flop (\ref{warn.flip-flop}).
\end{exrc}

\begin{exrc} Let $f:X\map Y$ be a birational map between
smooth, proper varieties. Show that
$$
\rho(X)-\rho(Y)=
\#\{\mbox{divisors contracted by $f$}\}-
\#\{\mbox{divisors contracted by $f^{-1}$}\}
$$
\end{exrc}

We are not yet ready to define minimal models.
As a first approximation, let us focus on
the codimension 1 part.

\begin{tempdef} Let $X$ be a smooth projective variety.
We say that $X$ is {\it minimal in codimension 1}
if every birational map $f:Y\map X$ 
from a smooth variety $Y$ is a birational  contraction.

In particular, this implies that $X$ has the smallest Picard number in its
birational equivalence class.
\end{tempdef}

\begin{exrc}  1. Let $X$ be a smooth projective variety
such that $K_X$ is nef. Prove that $X$ is minimal  in codimension 1.

2. $\p^3$ has the smallest Picard number in its
birational equivalence class but it is not  minimal  in codimension 1.

3. Let $X\subset \p^4$ be a smooth degree 4 hypersurface.
Then $K_X$ is not nef but, as proved by
Iskovskikh-Manin,  $X$ is  minimal  in codimension 1.
(See \cite[Chap.5]{ksc} for a proof and an introduction to 
these techniques.) 
\end{exrc}

\begin{exrc} Set
$X_0:=(x_1x_2+x_3x_4+x_5x_6=0)\subset  \a^6$.
Let $L\subset X_0$ be any 3-plane through the origin. 
Prove that, after a suitable coordinate change, 
$L$ can be given as $(x_1=x_3=x_5=0)$.
Prove that
$B_LX_0$ is smooth.

 Let
$Y$ be a smooth projective variety of dimension $3$
and $f_i, g_i$ are general sections of a very ample line bundle $L$ on $Y$.
Set 
$$
X':= \Bigl(\tsum_{i=1}^3 f_i({\mathbf x})g_i({\mathbf y})=0\Bigr)
\subset Y_{\mathbf x}\times Y_{\mathbf y},
$$
where ${\mathbf x}$ (resp.\ ${\mathbf y}$)
are the coordinates on the first (resp.\ second) factor.

  Assume that  $Y$
 admits a 
 finite morphism to an Abelian variety. Show that
$X'$ is  not birational to any smooth 
proper variety $X$
that is  minimal  in codimension 1.
\end{exrc}

\begin{exrc}[Contractions of products] \cite{ko-la}
 Let $X, U,V$ be normal projective varieties
and  $\phi:U\times V\map X$  a birational  contraction. 
 Assume that
$X$ is smooth (or at least has rational singularities).
Prove that  there are  normal projective varieties
$U'$ birational to $U$ and $V'$ birational to $V$ such that
$X\cong U'\times V'$.

In particular, $U\times V$ is minimal  in codimension 1
iff $U$ and $V$ are both  minimal  in codimension 1

\medskip

Hints to the proof. 
First reduce to the case when $U,V$  are smooth.

Let $|H|$ be a complete, very ample linear system on $X$ and $\phi^*|H|$
its pull back to $U\times V$.
Using that $\phi$ is a contraction, prove that
 $\phi^*|H|$ is also  a complete linear system.

 If $H^1(U,\o_U)=0$, then
$\pic(U\times V)=\pi_U^*\pic(U)+ \pi_V^*\pic(V)$, thus
there are divisors $H_U$ on $U$ and $H_V$ on $V$ such that
$\phi_*|H|\sim \pi_U^*H_U+\pi_V^*H_V$. Therefore 
$$
H^0(U\times V, \o_{U\times V}(\phi^*|H|))=
H^0(U,\o_U(H_U))\otimes H^0(V,\o_V(H_V)).
$$
Let now  $U'$ be the image of $U$
under the complete linear system $|H_U|$ and $V'$  the image of $V$
under the complete linear system $|H_V|$.

The $H^1(U,\o_U)\neq 0$ case is a bit harder.
Replace $H$ by a divisor
$H^*:=H+B$ where $B$ is a pull back of a divisor from the
product of the Albanese varieties of $U$ and $V$.
Show that for suitable $B$, 
there are divisors $H_U$ on $U$ and $H_V$ on $V$ such that
$\phi_*|H^*|\sim \pi_U^*H_U+\pi_V^*H_V$.
The rest of the argument now works as before.
\end{exrc}

\section{The cone of curves}

For details, see \cite[Chap.3]{km-book}.
\medskip

\begin{defn} Let $X$ be a projective variety over $\c$.
Any irreducible curve $C\subset X$ has a
homology class $[C]\in H_2(X,\r)$.
These classes generate a cone $NE(X)\subset H_2(X,\r)$,
called the {\it cone of curves} of $X$.
Its closure is denoted by $\nec(X)\subset H_2(X,\r)$.

If $X$ is over some other field, we can use
the vector space $N_1(X)$ of curves modulo numerical equivalence
instead of $H_2(X,\r)$ to define the
cone of curves $NE(X)\subset N_1(X)$.
\end{defn}

\begin{exrc} Show that every effective curve in
$\p^{a_1}\times\cdots\times \p^{a_n}$
is rationally equivalent to a nonnegative linear combination of
lines in the factors. Thus
$$
\nec\bigl(\p^{a_1}\times\cdots\times \p^{a_n}\bigr)\subset \r^n
$$
is the polyhedral cone spanned by the basis elements
corresponding to the lines.
\end{exrc}

\begin{exrc} 
Assume that a connected, solvable group acts on
$X$ with finitely many orbits. Show that
$\nec(X)$ is the polyhedral cone spanned by the
homology classes of the 
1-dimensional orbits. (The same holds even for rational equivalence
instead of homological equivalence.)

Hint. Use  the Borel fixed point theorem: 
A connected, solvable group acting on a proper variety has
a fixed point. Apply this to the Chow variety or the Hilbert
scheme parametrizing curves in $X$.
\end{exrc}

\begin{exrc} Let $S\subset \p^3$ be a smooth cubic surface.
Show that every effective curve is linearly equivalent to a linear
combination of lines. Thus
$\nec(X)\subset \r^7$
is a polyhedral cone spanned by the classes of the 27 lines.
(Note that the Cone theorem implies this only with rational
coefficients, not with integral ones. The proof is easiest using
the basic theory of linear systems.)
\end{exrc}

\begin{exrc}\label{EE.erc}
Let $A$ be an Abelian surface. 
If $Z$ is an ample $\r$-divisor, then $(Z\cdot Z)>0$.
Prove that, conversely, the condition
$(Z\cdot Z)>0$ defines a subset of $N_1(A)$  with 2 connected
components, one of which consists of ample $\r$-divisors.
Show that its closure is 
$\nec(A)$.

Check that if $A=E\times E$ and $E$ does not have complex multiplication
then every curve is algebraically equivalent 
to a linear combination $aE_1+bE_2+cD$ where $E_i$ are the two factors
and $D$ the diagonal. Thus
$$
\nec(E\times E)=\{aE_1+bE_2+cD: ab+bc+ca\geq 0 \mbox{ and } a+b+c\geq 0\}
\subset \r^3
$$
is a ``round'' cone.
\end{exrc}

Despite what these examples suggest, the cone of curves is
usually extremely difficult to determine.
For instance, we still don't know the
cone of curves for the following examples.
\begin{enumerate}
\item $C\times C$ for a general curve $C$.
(See \cite[Sec.1.5]{laz-book} for the known results and references.)
\item The blow up of $\p^n$ at more than a few points, cf.\
\cite{cast}.
\end{enumerate}

A basic discovery of \cite{mori-tf} is that
the part of the cone of curves which has negative intersection
with the canonical class is quite well behaved.
Subsequently it was generalized to certain perturbations of
the canonical class. The precise definitions
will be given in Section \ref{sing.sec}. For now you can
imagine that $X$ is smooth and
$\Delta=\sum a_i D_i$ is a $\q$-divisor
where $\sum D_i$ is a simple normal crossing divisor
and $0<a_i<1$ for every $i$.

\begin{thm}[Cone theorem] (cf.\ \cite[Thm.3.7.1--2]{km-book}) \label{cone1.thm}
	Let $(X,\Delta)$ be a projective klt pair with $\Delta$
effective.  Then: 
\begin{enumerate}
\item There are (at most countably many) rational curves
$C_j\subset X$ such that  
$0<-(K_X+\Delta)\cdot C_j\leq 2\dim X$ and
$$
\nec(X) =
\nec(X)_{(K_X+\Delta)\geq 0} + \sum
\r_{\geq 0}[C_j],
$$
where $\nec(X)_{(K_X+\Delta)\geq 0}$ denotes the set of those
elements of $\nec(X)$ that have nonnegative
intersection number with $K_X+\Delta$.

\item  For any $\epsilon > 0$ and ample $\q$-divisor $H$, 
	$$
\nec(X)  
	           =  \nec(X)_{(K_X+\Delta+\epsilon H)\geq
0}  +
\sum_{{\rm finite}}\r_{\geq 0}[C_j].
	$$
\end{enumerate}
\end{thm}

If $-(K_X+\Delta)$ is ample then taking
$H=-(K_X+\Delta)$ and $\epsilon<1$ in (\ref{cone1.thm}.2),
the first summand on the right is trivial.
Hence we obtain:

\begin{cor} Let $(X,\Delta)$ be a projective klt pair with $\Delta$
effective and $-(K_X+\Delta)$  ample. 
There are finitely many rational curves
$C_j\subset X$ such that  
	$$
\nec(X)  =\sum\r_{\geq 0}[C_j].
	$$
In particular, $\nec(X)$ is a polyhedral cone.\qed
\end{cor}

\begin{warning} If the cone is 3-dimensional, the cone theorem
implies that the $(K_X+\Delta)$-negative part of 
$\nec(X)$ is locally polyhedral. This, however, fails 
for 4-dimensional cones.

Use (\ref{EE.erc}) to show that such an example is given by
$\nec(E\times E\times\p^1)$
where $E$ is an elliptic curve which does not have complex multiplication.
\end{warning}

\begin{defn}\label{exray.defn} In convex geometry, a closed subcone
$F\subset \nec(X)$ is  called an {\it extremal face}
if $u,v \in  \nec(X)$ and $u+v\in F$ implies that $u,v\in F$.
A 1-dimensional extremal face is called an {\it extremal ray}.

In algebraic geometry, one frequently assumes in addition
that intersection product with 
$K_X$ (or $K_X+\Delta$) gives a  strictly negative linear function
on $F\setminus\{0\}$.

Thus, extremal rays of $\nec(X)$ are precisely those
 summands $\r_{\geq 0}[C_j]$ in (\ref{cone1.thm}.1)
that are actually
needed. 
\end{defn}

The next result shows that
there are contraction morphisms associated to
any extremal face.

\begin{thm}[Contraction theorem] (cf.\ \cite[Thm.3.7.2--4]{km-book})
\label{cone2.thm}
	Let $(X,\Delta)$ be a projective klt pair with $\Delta$
effective.  Let
$F\subset \nec(X)$ be a ($(K_X+\Delta)$-negative)
extremal face. Then there is a unique   morphism
$\cont_F:X\to Z$, called the {\it contraction}  of $F$,   such that
$(\cont_F)_*\o_X=\o_Z$ and  an irreducible curve $C\subset X$
is mapped to  a point by $\cont_F$ iff
$[C]\in F$.
Moreover, 
\begin{enumerate}
\item $R^i(\cont_F)_*\o_X=0$ for $i>0$, and
\item if  $L$ is a
line bundle on $X$ such that $(L\cdot C)=0$ whenever $[C]\in F$ 
 then there is a line bundle $L_Z$ on
$Z$ such that $L\cong \cont_F^*L_Z$. 
\end{enumerate}
\end{thm}

\begin{exrc} Let $Z$ be a smooth, projective variety and
$W\subset X$ a smooth, irreducible subvariety of codimension $\geq 2$.
Show that $\pi:B_WZ\to Z$ is the 
contraction of an extremal ray on $B_WZ$.
\end{exrc}

\begin{exrc} Let $Z$ be an $n$-dimensional   projective variety 
with a unique singular point $p$ of the form
$$
x_1^m+\cdots+x_{n+1}^m+(\mbox{higher terms})=0.
$$
Show that $B_pZ$ is smooth and $\pi:B_pZ\to Z$ is the 
contraction of an extremal face on $B_pZ$
iff $m<n$. The exceptional divisor
is the smooth hypersurface
$(x_1^m+\cdots+x_{n+1}^m=0)\subset \p^{n}$.

If $n\geq 4$, then by the Lefschetz theorem, 
$\pi$ is the 
contraction of an extremal ray. 
Find examples with $n=3$ where we do contract a face.
\end{exrc}

\begin{exrc}\label{rho=1.gen.exmp} Let $f_i(x_1,\dots,x_4)$ for $i=m,m+1$ be
homogeneous of degree $i$. Assume that
$$
X:=\bigl(x_0f_m(x_1,\dots,x_4)+f_{m+1}(x_1,\dots,x_4)=0\bigr)\subset \p^4
$$
is smooth away from the origin.
Prove that every Weil divisor on $X$ is obtained by
intersecting $X$ with another hypersurface.
\end{exrc}

\begin{exrc} Let $Z$ be an $n$-dimensional   projective variety 
with a unique singular point $p$ of the form
$$
x_1^m+\cdots+x_n^m+x_{n+1}^{m+1}+(\mbox{higher terms})=0.
$$
Show that  $B_pZ$ is smooth and $\pi:B_pZ\to Z$ is the 
contraction of an extremal ray on $B_pZ$
iff $m<n$ and $n\geq 3$. The exceptional divisor
is the singular  hypersurface
$(x_1^m+\cdots+x_{n}^m=0)\subset \p^{n}$.
\end{exrc}

\begin{exrc} Let $f_m(x_1,\dots,x_{n+1})$ be an irreducible,
homogeneous degree $m$ polynomial and 
$g_{m+1}(x_1,\dots,x_{n+1})$  a general,
homogeneous degree $m+1$ polynomial.
Let $Z$ be an $n$-dimensional   projective variety 
with a unique singular point $p$ of the form
$$
f_m(x_1,\dots,x_{n+1})+g_{m+1}(x_1,\dots,x_{n+1})+(\mbox{higher terms})=0.
$$
Use (\ref{bert.exrc}) and (\ref{cA.can.exrc}) 
to prove that $B_pZ$ has only
canonical singularities (\ref{7.9}).

Show that $\pi:B_pZ\to Z$ is the 
contraction of an extremal face on $B_pZ$
iff $m<n$. 

Note that the exceptional divisor
is the   hypersurface
$(f_m(x_1,\dots,x_{n+1})=0)\subset \p^{n}$,
which can be quite singular.
\end{exrc}

\begin{exrc} Let $Z\subset \p^n$ be defined by
$x_0=f(x_1,\dots,x_n)=0$ where $f$ is irreducible.
Show that $B_Z\p^n\to\p^n$ is the 
contraction of an extremal ray on $B_Z\p^n$.
Show that $Z$ has only $cA$-type singularities (\ref{cA.can.exrc}).
When is $Z$ canonical or terminal (\ref{7.9})?

Note that the exceptional divisor
is a $\p^1$-bundle over $Z$,
which can be quite singular.
\end{exrc}

\begin{exrc}\label{gen.cast.erc}  Let $X$ be a smooth, projective variety,
$D\subset X$ a smooth hypersurface and
$C\subset D$ any curve.
Assume that the Picard number of $D$ is 1 and
the conormal bundle  $N^*_{D|X}$ is ample.

Prove that $[C]$ is an  extremal ray of $\nec(X)$
in the convex geometry sense (\ref{exray.defn}).
When is it a $K_X$-negative extremal ray?

Assume in addition that $-K_D$ is ample.
 Generalize the proof of Castelnuovo's theorem
(for instance, as in \cite[V.5.7]{hartsh})
to prove  (\ref{cone2.thm})  in this case. (That is, 
 there is a contraction
$\pi:X\to X'$ that maps $D$ to a point and is an isomorphism on
$X\setminus D$.)
\end{exrc}

\begin{exrc} With notation as in (\ref{gen.cast.erc}), assume
 that  $D\cong \p^{n-1}$
and $N_{D|X}\cong \o_{\p^{n-1}}(-m)$.
Set $x':=\pi(D)$. Prove that the completion of
$X'$ (at $x'$)
  is  isomorphic to the completion (at the origin) of the quotient 
of $\a^n$ by the $\z/m$-action
$(x_1,\dots,x_n)\mapsto (\epsilon x_1,\dots,\epsilon x_n)$
where $\epsilon $ is a primitive $m$th root of 1.
(Hint: Use the methods of \cite[Exrc.II.8.6--7]{hartsh}.)
\end{exrc}

\begin{exrc} Let $Z$ be a smooth, projective variety and
$X\subset Z\times \p^m$ a smooth hypersurface
such that $X\cap\bigl(\{z\}\times \p^m\bigr)$
is a hypersurface of degree $d$ for general $z\in Z$.

Show that the projection $\pi:X\to Z$ is the 
contraction of an extremal face on $X$
iff $d<m+1$ and $m\geq 2$.
 
If $m=2$ and $\dim Z=2$ then show that
every fiber of $\pi:X\to Z$ is either a line
(if $d=1$) or a (possibly singular) conic
(if $d=2$). (This can fail if $X$ has an ordinary double point.)

If  $m=2$ and $\dim Z=3$ then 
find smooth examples where the general 
 fiber of $\pi:X\to Z$ is  a line or a 
 conic but special fibers are $\p^2$.
\end{exrc}

\begin{exrc}\label{1-dim.fiber.contr} If you know some about the
deformation theory and the Hilbert scheme 
 of curves on smooth varieties, prove the
following.
 (You will find (\ref{cone2.thm}.1)  very helpful.)

Let $\pi:X\to Z$ be an extremal contraction
with $X$ smooth
where every fiber has dimension $\leq 1$. 
Then $Z$ is smooth and we have one of the following cases:
\begin{enumerate}
\item $X=B_WZ$ for some smooth
$W\subset Z$ of codimension 2.
\item $X$ is a $\p^1$-bundle over $Z$.
\item $X$ is a hypersurface in a $\p^2$-bundle over
$Z$ and every fiber of $\pi:X\to Z$ is  a (possibly singular) conic.
\end{enumerate}
\end{exrc}

\begin{exrc} Let $X\subset \p^4$ be a degree 3 hypersurface with a
unique singular point that is an ordinary node.
(That is, analytically isomorphic to $(xy-zt=0)$.)

Let $\pi:Y\to X$ denote the blow up of the node.
Prove that its exceptional divisor $E$ is isomorphic to
$\p^1\times \p^1$ and its normal bundle is $\o_{\p^1\times \p^1}(-1,-1)$.

Thus $E$ looks like it could have been obtained
by blowing up a curve $C\cong \p^1$
with normal bundle $\o_{\p^1}(-1)+\o_{\p^1}(-1)$ in a smooth 3-fold.
Nonetheless, use (\ref{rho=1.gen.exmp}) to 
show that there is no such projective 3-fold.
\end{exrc}

\begin{exmp}  Let $X$ be the $cE_7$-type singularity
$(x^2+y^3+yg_3(z,t)+h_5(z,t)=0)\subset \a^4$,
 where $g_3$ and $h_5$ do not have a
common factor. Show  that
$X$ has an isolated singular point at the origin and its
$(3,2,1,1)$-blow up $Y\to X$ has only terminal singularities. 
(See \cite[4.56]{km-book} or \cite[6.38]{ksc} for
weighted blow-ups.) 
Conclude from this that $X$ itself has a terminal singularity.

One of the standard  charts on the blow up is given by the
substitutions
$x=x_1y_1^3, y=y_1^2, z=z_1y_1, t=t_1y_1$ and the
exceptional divisor has equation
$$
 E=(g_3(z_1,t_1)+h_5(z_1,t_1)=0)/
\tfrac12(1,1,1) \subset \a^3/\tfrac12(1,1,1).
$$ 
This gives  examples of extremal contractions whose exceptional
divisor $E$ has  quite complicated singularities.
\begin{enumerate}
\item $x^2+y^3+yz^3+t^5$. $E$ is singular along 
$(z_1=t_1=0)$, with a transversal singularity type
$z^3+t^5$, that is $E_8$.
\item $x^2+y^3+y(z-at)(z-bt)(z-ct)+t^5$. $E$ has  triple
self-intersection  along 
$z_1=t_1=0$.
\end{enumerate}
\end{exmp}

\begin{exrc}  Let $X$ be a smooth Fano variety, $\dim X\geq 4$. Let
$Y\subset X$ be a smooth divisor in $|-K_X|$
(thus $K_Y=0$).
Show that the natural map $ i_*:\nec(Y)\to \nec(X) $
is an isomorphism. Thus $\nec(Y)$
is a polyhedral cone.
(See \cite{bor1, bor2} for many such interesting examples.)
\medskip

Steps of the proof.

1.  By a theorem of Lefschetz,
 $i_*$ is an injection. Thus we need to show that
for every extremal ray $R$ of $\nec(X)$ there is a curve
 $C_R\subset Y$ such that
$C_R$ generates $R$ in $\nec(X)$.

2. Let $f:X\to Z$ be the contraction morphism of $R$. If there is a fiber
$F\subset X$ of $f$ whose dimension is at least two then $Y\cap F$ contains a
curve $C_R$ which works.

3. If every fiber of $f$ has dimension one then we use 
(\ref{1-dim.fiber.contr}).
We need to show that in these cases $Y$ contains a fiber of $f$.

4. Prove the following lemma.
 Let $g:U\to V$ be a $\p^1$-bundle over a normal
projective variety.  Let $V'\subset U$ be an irreducible divisor
 such that $g:V'\to V$ is
finite of degree one (thus an isomorphism). If $V'$ is ample then
$\dim V\leq 1$.

5. 
In the divisorial contraction case 
apply this lemma to  $U:=$ the exceptional divisor of $f$. 

6. In the $\p^1$-bundle case
apply this lemma to $U:=$ normalization of the branch divisor of $Y\to Z$.
(If there is no branch divisor, then to $X\times_ZY\to Y$.)

7. In the conic bundle case there are two possibilities.
If every fiber is smooth, this is like the $\p^1$-bundle case.
Otherwise apply the lemma to $U:=$ normalization of the divisor of 
singular fibers of $Y\to Z$.
\end{exrc}

\begin{exrc}\label{bert.exrc}
Prove the following result of \cite[4.4]{k-pairs}.

 {\it Theorem.} Let $X$ be a smooth variety over a field of
characteristic zero and  $|B|$  
a linear system of Cartier divisors.
Assume that for every $p\in X$ there is a $B(p)\in |B|$ such that $B(p)$ is
 smooth at $p$ (or $p\not\in B(p)$). 

Then a   general  member 
$B^g\in |B|$ has only  $cA$-type singularities (\ref{cA.can.exrc}).

Hint.
 By Noetherian induction it is sufficient to prove
that for every irreducible subvariety $Z\subset X$ there is an open subset
$Z^0\subset Z$ such that a   general  member 
$B^g\in |B|$ has only  $cA$-type singularities at points of $Z^0$.

If $Z\not\subset \bs|B|$ then use the usual Bertini
theorem.

If  $Z\subset\bs|B|$ and $\codim(Z,X)=1$, then use the usual Bertini
theorem for $|B|-Z$.

If  $Z\subset\bs|B|$ and $\codim(Z,X)>1$ then restrict to a
suitable hypersurface
 $Z\subset Y\subset X$ and use induction.
\end{exrc}

\begin{exrc} 
Use the following examples to show that the conclusion of (\ref{bert.exrc})
is almost optimal:

 Let $X=\c^n$ and $f\in \c[x_3,\dots,x_n]$
such that    $(f=0)$ has an isolated singularity at
the origin. Consider the linear system $|B|=(\lambda x_1+\mu x_1x_2+
\nu f=0)$. Show that at each point there is a smooth
member and the general member is
singular at 
$(0,-\lambda/\mu,0,\dots,0)$ with local equation
$(x_1x_2+f=0)$. 

  Consider the linear system $\lambda(x^2+zy^2)+\mu y^2$. At any point
$x\in \c^3$ its general member has a $cA$-type singularity, but the general
member has a moving pinch point. 
\end{exrc}

\section{Singularities}\label{sing.sec}

For details, see \cite[Chaps.4--5]{km-book}. 
\medskip

We already saw in several examples that
even if we start with a smooth variety, the contraction
of an extremal ray can lead to a singular variety.
It took about 10 years to understand the correct classes
of singularities that one needs to consider.
Instead of going through this historical process,
let us jump into the final definitions.

\begin{rem} In the early days of the MMP, a lot of effort was devoted
to classifying the occurring singularities in  dimensions 2 and 3.
While it is comforting to have some concrete examples and
lists at hand, the recent advances use very little of these
explicit descriptions. In most applications, we fall back to the
definitions via log resolutions. The key seems to be an ability
to work with  log resolutions.
\end{rem}

\begin{defn}\label{mmp.discr.def} Let $X$ be a normal scheme
and $\Delta$ a $\q$-divisor on $X$ such that
$K_X+\Delta$ is $\q$-Cartier. Let
$f:Y\to X$  be a birational morphism, $Y$ normal.
Let $E_i\subset \ex(f)$ be the exceptional divisors. If
$m(K_X+\Delta)$ is Cartier, then 
$f^*\o_X\bigl(m(K_X+\Delta)\bigr)$ is defined and there is a natural isomorphism
$$ f^*\o_X\bigl(m(K_X+\Delta)\bigr)|_{Y\setminus \ex(f)}\cong 
\o_Y\bigl(m(K_Y+f^{-1}_*\Delta)\bigr)|_{Y\setminus \ex(f)},
\eqno{(\ref{mmp.discr.def}.1)}
$$
where $f^{-1}_*\Delta$ denotes the birational transform of $\Delta$.
 Hence there are integers $b_i$ such that
$$
\o_Y\bigl(m(K_Y+f^{-1}_*\Delta)\bigr)\cong 
f^*\o_X\bigl(m(K_X+\Delta)\bigr)(\tsum b_iE_i).
\eqno{(\ref{mmp.discr.def}.2)}
$$ 
Formally divide by $m$ and write this as
$$
 K_Y+\Delta_Y\simq f^*(K_X+\Delta)
\qtq{where} \Delta_Y:=f^{-1}_*\Delta-\tsum (b_i/{m})E_i.
$$
 The rational number $a(E_i,X,\Delta):=b_i/m$ is called the {\it discrepancy} of
$E_i$ with respect to $(X,\Delta)$.

The closure of $f(E_i)\subset X$ is called the {\it center} of $E_i$ on
$X$. It is denoted by $\cent_XE_i$.

If $f':Y'\to X$ is another birational morphism and 
$E'_i:=\bigl((f')^{-1}\circ f\bigr)(E_i)\subset Y'$ 
is a divisor then $a(E'_i,X,\Delta)=a(E_i,X,\Delta)$
 and
$\cent_XE_i=\cent_XE'_i$.  Thus the discrepancy and the center depend
only on the divisor up to birational equivalence, but not on the
particular variety where
 the divisor appears.
\end{defn}

\begin{defn}  Let $X$ be a normal variety. An {\it $\r$-divisor}
on $X$ is a formal $\r$-linear combination $\sum r_iD_i$
of Weil divisors. We say that two $\r$-divisors $A_1, A_2$ are
{\it $\r$-linearly equivalent,} denoted $A_1\simr A_2$,
 if  there are rational functions
$f_i$ and real numbers $r_i$ such that
 $A_1-A_2=\sum r_i(f_i)$.

One can pretty much work with $\r$-divisors as with
$\q$-divisors, but some basic properties need to be
thought through.
\end{defn}

\begin{exrc} Prove the following about  $\r$-divisors and
$\r$-linear equivalence.

(1) Let $A_1, A_2$ be two $\q$-divisors. Show that
$A_1\simr A_2$ iff $A_1\simq A_2$.

(2) Define the pull back of $\r$-divisors and show that it is well defined.

(3) Let $A$ be an $\r$-divisor such that $A\simr 0$.
Prove that one can write $A=\sum r_i(f_i)$ such that
$\supp \bigl((f_i)\bigr)\subset \supp A$ for every $i$.
\end{exrc}

\begin{exrc} 
Let $X$ be a normal scheme
and $\Delta$ an $\r$-divisor on $X$ such that
$K_X+\Delta$ is $\r$-Cartier. Let
$f:Y\to X$  be a proper birational morphism, $Y$ normal.
 Show that there is a  unique $\r$-divisor 
$\Delta_Y$ such that
\begin{enumerate}
\item $f_*\bigl(\Delta_Y\bigr)=\Delta$, and
\item $ K_Y+\Delta_Y\equiv_f f^*(K_X+\Delta)$, where
$\equiv_f$ denotes relative numerical equivalence, that is,
 $ (K_Y+\Delta_Y\cdot C)=( f^*(K_X+\Delta)\cdot C)$ for every curve
 $C\subset Y$
such that $\dim\ f(C)=0$. (Note that the latter is just $0$.)
\end{enumerate}
Use this to define discrepancies for $\r$-divisors.
\end{exrc}

\begin{exrc} Formulate (\ref{mmp.discr.def})
in case $f:Y\map X$ is a birational map which is defined
outside a codimension 2 set. (This holds, for instance
if $X$ is proper over the base scheme $S$.)
\end{exrc}

\begin{exrc}[Divisors and rational maps]
Let $f:X\map Y$ be a generically finite rational map between
proper, normal schemes.
Define the push forward $f_*:\wdiv(X)\to \wdiv(Y)$
of Weil divisors. Show that if $f,g$ are
morphisms then $(f\circ g)_*=f_*\circ g_*$ but this fails
even for birational maps.

Let $f:X\map Y$ be a  dominant rational map between
normal schemes, $Y$ proper.
Define the pull back $f^*: \cdiv(Y)\to \wdiv(X)$
from Cartier divisors to Weil divisors.
Show that if $f$ is a morphism then we get
$f^*: \cdiv(Y)\to \cdiv(X)$ but not in general.
Find examples of birational maps between smooth projective varieties
such that $(f\circ g)^*\neq f^*\circ g^*$.
\end{exrc}

\begin{defn}\label{7.9}Let $(X,\Delta)$ be a pair where $X$ is a
normal variety and
$\Delta=\sum a_iD_i$ is a sum of distinct prime divisors. (We allow
the $a_i$ to be arbitrary real numbers.) Assume that 
 $K_X+\Delta$ is $\r$-Cartier.  We say that
$(X,\Delta)$ is
$$
\left.
\begin{array}{c}
 \mbox{\it terminal}   \\       
 \mbox{\it canonical}  \\
   \mbox{\it klt}  \\
   \mbox{\it plt}  \\
\begin{array}{l}
{\ }\\
{\ } 
\end{array}
\mbox{\it dlt}
\begin{array}{l}
{\ }\\
{\ } 
\end{array}\\
\mbox{\it lc}    
\end{array}
\right\}
 \quad\mbox{if $a(E,X,\Delta)$ is}\quad 
\left\{
\begin{array}{l}
 >0\quad \forall\ E \mbox{ exceptional,}  \\       
 \geq 0\quad \forall\ E \mbox{ exceptional,}  \\
  >-1 \quad \forall\ E, \\
  > -1\quad \forall\ E \mbox{ exceptional,}  \\
  >-1 \
\begin{array}{l}
\forall\ E  \mbox{ such that $(X,\Delta)$ is not snc at}\\
\quad\mbox{the generic point of $\cent_X(E)$,}
\end{array}\\
\geq -1  \quad \forall\ E.
\end{array}
\right.                   
$$\end{defn}
 
Here klt is short for 
{\it Kawamata log terminal}, plt for  {\it purely log terminal},
dlt for  {\it divisorial log terminal},
 lc for {\it log canonical} and snc for {\it simple normal crossing}.
(The  phrase      ``$(X,\Delta)$ has terminal etc.
singularities''  may be confusing since it could refer to the
singularities of
$(X,0)$ instead.)

Each of these 5 notions has an important place in the theory of
minimal models: 

\begin{enumerate}
\item {\it Terminal}.  Assuming $\Delta=0$, this is the smallest class that
is necessary to run the minimal model program for smooth varieties. 
If $(X,0)$ is terminal then $\sing X$ has codimension $\geq 3$.
All 3-dimensional terminal singularities are classified,
see (\ref{rt-crit}) for some examples.
It is generally believed that already in dimension 4
a complete classification would be impossibly long.
The $\Delta\neq 0$ case appears only infrequently.

\item  {\it Canonical}.  Assuming $\Delta=0$, these are precisely the
singularities that appear on the canonical models of varieties of
general type. 
Two dimensional canonical singularities are classified,
see (\ref{dv.list}). There is some structural information in dimension
3 \cite[5.3]{km-book}. 
This class  is especially important for moduli problems.

\item {\it Kawamata log terminal}.  This is the smallest class that
is necessary to run the minimal model program for pairs
 $(X,\Delta)$ where $X$ is  smooth and $\Delta$ a
simple normal crossing divisor with coefficients $< 1$.

The vanishing theorems (cf.\ \cite[2.4--5]{km-book})
seem to hold naturally in this class. In general,  proofs that work
with canonical singularities frequently work with klt.
Most unfortunately, this class is
 not large enough  for inductive proofs.  

\item{\it Purely log terminal}. This is useful mostly for
inductive purposes. $(X,\Delta)$ is plt iff 
$(X,\Delta)$ is dlt and the irreducible components of $\rdown{\Delta}$ 
are disjoint. 

\item {\it Divisorial     log terminal}. This is the smallest class that
is necessary to run the minimal model program for pairs
 $(X,\Delta)$ where $X$ is  smooth and $\Delta$ a
simple normal crossing divisor with coefficients $\leq 1$.

By \cite{sza}, there is a
log resolution $f:(X',\Delta')\to (X, \Delta)$ such that
every $f$-exceptional divisor has discrepancy $>-1$
and $f$ is an isomorphism over the snc locus of
$(X, \Delta)$.

While the definition of this class is somewhat artificial looking,
it has good cohomological properties and is much better behaved
than general log canonical pairs.

If $\Delta=0$ then the notions klt and dlt 
coincide and in this case we say that $X$ has 
{\it log terminal} singularities
(abbreviated as {\it lt}).

\item {\it Log canonical}.  This is the largest class where discrepancy
still makes sense and inductive arguments naturally run in this class.
There are three major complications though:
\begin{enumerate}
\item The refined vanishing theorems  frequently fail.
\item The singularities are not rational and not even Cohen-Macaulay,
hence rather complicated
from the cohomological point of view; see, for example, (\ref{rt-crit}).
\item Several tricks of perturbing coefficients can not be done
since a perturbation would go above $1$;
see, for example, (\ref{nonfg.exrc}). 
\end{enumerate}
\end{enumerate}

\begin{exrc} \label{discr.up.exrc} Let $f:X\to Y$ be 
a birational morphism,  $\Delta_X$, $\Delta_Y$ $\r$-divisors 
such that $f_*\Delta_X=\Delta_Y$
and $D$ 
an  effective $\r$-divisor. Assume that 
$K_Y+\Delta_Y$ and $D$ are $\r$-Cartier and
$$
K_X+\Delta_X\simr f^*(K_Y+\Delta_Y)+ D.
$$
Prove that for any $E$,
$a(E,X,\Delta_X)\leq a(E,Y,\Delta_Y)$
and the inequality is strict iff $\cent_XE\subset \supp D$.
\end{exrc}

\begin{exrc}
Show that the assumptions of (\ref{discr.up.exrc}) are fulfilled 
(for suitable $\Delta_Y$ and $D$) if
$X$ is $\q$-factorial, 
$f$ is 
the birational contraction of a $(K_X+\Delta_X)$-negative extremal ray
and $\ex(f)$ has codimension 1.
\end{exrc}

The following exercise shows why
log canonical is the largest class defined.

\begin{exrc} Given $(X,\Delta)$ assume that there is
a divisor $E_0$ such that $a(E_0,X,\Delta)<-1$.
Prove that  $\inf_E\{a(E,X,\Delta)\}=-\infty$.
\end{exrc}

\begin{exrc}
Show that if $(X,\sum a_iD_i)$ is lc  (and the $D_i$ are distinct) then
$a_i\leq 1$ for every $i$.
\end{exrc}

\begin{exrc} Assume that $X$ is smooth and $\Delta$ is effective.
Show that if   $\mult_x\Delta <1$ (resp.\ $\leq 1$)
for every $x\in X$ then $(X,\Delta)$ is terminal (resp.\ canonical).

Prove that the converse holds for surfaces but not in higher dimensions.
\end{exrc}

\begin{exrc}[Du Val singularities]\label{dv.list}

In each of the following cases, construct the minimal resolution and
verify that its dual graph is the graph given.
Check that these singularities are canonical.
(One can see that these are all the
2-dimensional canonical singularities.)
 See \cite[Sec.4.2]{km-book} or \cite{durfee}
for more information.
(The equations below are correct in characteristic zero.
The dual graphs are correct in every characteristic.)
\medskip

\noindent$A_n$: $x^2+y^2+z^{n+1}=0$, with $n\geq 1$ curves in the dual graph:
$$
\begin{array}{ccccccccc}
2 & - & 2 & - & \cdots &- & 2 & - & 2
\end{array}
$$
$D_n$: $x^2+y^2z+z^{n-1}=0$, with $n\geq 4$ curves in the dual graph:
$$
\begin{array}{ccccccccc}
&& 2&&&&&&\\
&& \vert &&&&&&\\
2 & - & 2 & - & \cdots &- & 2 & - & 2
\end{array}
$$
$E_6$: $x^2+y^3+z^4=0$, with   dual graph:
$$
\begin{array}{ccccccccc}
&& &&2&&&&\\
&&&& \vert &&&&\\
2 & - & 2 & - & 2 &- & 2 & - & 2
\end{array}
$$
$E_7$: $x^2+y^3+yz^3=0$, with  dual graph:
$$
\begin{array}{ccccccccccc}
&& &&2&&&&&&\\
&&&& \vert &&&&&&\\
2 & - & 2 & - & 2 &- & 2 & - & 2& - & 2
\end{array}
$$
$E_8$: $x^2+y^3+z^5=0$, with  dual graph:
$$
\begin{array}{ccccccccccccc}
&& &&2&&&&&&&&\\
&&&& \vert &&&&&&&&\\
2 & - & 2 & - & 2 &- & 2 & - & 2& - & 2& - & 2
\end{array}
$$
\end{exrc}

\begin{exrc}[$cA$-type singularities]\label{cA.can.exrc}
Let $0\in X$ a normal $cA$-type singularity.
That is, 
 either $X$ is
 smooth at $0$, or, in suitable local coordinates $x_1,\dots, x_n$,
the equation of $X$ is $x_1x_2+(\mbox{other terms})=0$.

Show that $X$ is 
\begin{enumerate}
\item canonical near $0$ iff $\dim \sing X\leq \dim X-2$, and
\item terminal near $0$ iff $\dim \sing X\leq \dim X-3$.
\end{enumerate}

Hint. First show that being $cA$-type is an open
condition. Then use a lemma of Zariski and Abhyankar
(cf.\ \cite[2.45]{km-book}) to reduce everything to the
statements:
\begin{enumerate}\setcounter{enumi}{2}
\item The exceptional divisor(s) of $B_0X\to X$ have
discrepancy $\dim X-2$, save when $X$ is smooth.
\item $B_0X$ has only $cA$-type singularities.
\end{enumerate}
\end{exrc}

\begin{exrc}[Some simple elliptic singularities]\label{ell.lowmult.list}

In each of the following cases, construct the minimal resolution.
Verify that the exceptional set  is a single elliptic curve with
self intersection $-k$.
\begin{enumerate}
\item[($k=3$)]\quad   $(x^3+y^3+z^3=0)$.  (This is very easy)
\item[($k=2$)]\quad   $(x^2+y^4+z^4=0)$.
\item[($k=1$)]\quad  $(x^2+y^3+z^6=0)$.  (This is a bit tricky.)
\end{enumerate}

In general, prove that for any elliptic curve $E$ and any $k\geq 1$
there is a  normal singularity whose minimal resolution
contains $E$ as the  single exceptional curve with
self intersection $-k$.

Check that all of these are log canonical.

Use the methods of \cite[Exrc.II.8.6--7]{hartsh}
to prove that the completion of the singularity  is uniquely 
determined by $E$.
\end{exrc}

\begin{exrc}\label{ell.lowmult.quot.list}
 Construct the minimal resolutions of the following
quotients of the singularities in (\ref{ell.lowmult.list}).
(See (\ref{quot.sing.say}) for the notation.)
\begin{enumerate}
\item[] $(x^3+y^3+z^3=0)$:  $\tfrac13(1,0,0)$,  $\tfrac13(1,1,1)$.
\item[]   $(x^2+y^4+z^4=0)$: $\tfrac12(1,0,0)$,  $\tfrac14(0,0,1)$.
\item[]  $(x^2+y^3+z^6=0)$: $\tfrac16(0,0,1)$.
\end{enumerate}
\end{exrc}

\begin{exrc}\label{cone.gen.erc}
Let $X\subset \p^n$ be a smooth variety
and $C(X)\subset \a^{n+1}$ the cone over $X$.
Show that $C(X)$ is normal iff
$H^0(\p^n, \o_{\p^n}(m))\to H^0(X, \o_X(m))$ is onto
for every $m\geq 0$.

Assume next  that $C(X)$ is normal. 
Let $\Delta$ be an effective $\q$-divisor on $X$.
Prove that
\begin{enumerate}
\item $K_{C(X)}+C(\Delta)$ is $\q$-Cartier iff
$K_X+\Delta\simq r\cdot H$ for some $r\in\q$ where $H\subset X$ is the
hyperplane class.
\item If $K_X+\Delta\simq r\cdot H$ then  
$\bigl(C(X), C(\Delta)\bigr)$
 is
\begin{enumerate}
\item terminal iff $r<-1$ and $(X,\Delta)$ is terminal,
\item canonical iff $r\leq -1$ and $(X,\Delta)$ is canonical,
\item klt iff $r<0$ and $(X,\Delta)$ is klt, and
\item lc iff $r\leq 0$ and $(X,\Delta)$ is lc.
\end{enumerate}
\end{enumerate}
\end{exrc}

\begin{exrc}\label{rt-crit} Notation as in (\ref{cone.gen.erc}).
Prove that $C(X)$ has a rational singularity
iff $H^i(X,\o_X(m))=0$ for every $i>0, m\geq 0$
and a Cohen-Macaulay singularity iff
$H^i(X,\o_X(m))=0$ for every $\dim X>i>0, m\geq 0$.
In particular:
\begin{enumerate}
\item If $X$ is an Abelian variety and $\dim X\geq 2$
then $C(X)$ is log canonical but not Cohen-Macaulay.
\item If $X$ is a K3 surface
then $C(X)$ is log canonical, Cohen-Macaulay but not rational.
\item If $X$ is an Enriques surface
then $C(X)$ is log canonical and rational.
\end{enumerate}
\end{exrc}

\begin{say}[Quotient singularities]\label{quot.sing.say}
Let $G$ be any finite group. A homomorphism
 $G\to GL_n$ is equivalent to a
linear $G$-action  on $\a^n$.
The resulting quotient singularities $\a^n/G$ are rather special
but they provide a very good test class for many questions
involving log-terminal singularities.

One can always reduce to the case when
the $G$-action on $\a^n$ is effective and fixed point free outside 
a codimension 2 set. (Unless you are into stacks.)
Thus assume this in the sequel.

Show that any such $\a^n/G$ is log terminal.

Show that if $G\subset SL_n$ then
the canonical class of  $\a^n/G$ is  Cartier.
In particular,  $\a^n/G$  is canonical.

Assume that $G=\langle g\rangle$ is a cyclic group.
Any cyclic action on $\a^n$ can be diagonalized and written as
$$
g:(x_1,\dots, x_n)\mapsto (\epsilon^{a_1}x_1,\dots,  \epsilon^{a_n}x_n),
$$
where $\epsilon= e^{2\pi i/m}$, $m=|G|$  and $0\leq a_j<m$.
Define the {\it age} of $g$ as
$\operatorname{age}(g):=\frac1{m}(a_1+\cdots+a_n)$.
As a common shorthand notation, we denote the quotient by this action by
$$
\a^n/\tfrac1{m}(a_1,\dots, a_n).
$$
The following  very useful  criterion tells us
when  $\a^n/G$ is  terminal or canonical.
\medskip

{\it Reid-Tai  criterion.} $\a^n/G$  is canonical
(resp.\ terminal) iff the age of every non-identity element
$g\in G$ is $\geq 1$ (resp.\ $>1$).
\medskip

(This is not hard to prove if you know some basic toric techniques.
Otherwise, look up \cite{reid-yp}.)

As a consequence, prove that
the  3-fold quotients 
$\a^3/\tfrac1{m}(1,-1,a)$
are terminal
if  $(a,n)=1$. 
(It is a quite tricky combinatorial argument to show that
these are all the 3-dimensional terminal quotients, cf.\ \cite{reid-yp}.)

By contrast, every ``complicated'' higher dimensional
quotient singularity is terminal. By the results of
\cite{ko-la, gu-ti2}, if the $G$-action on $\a^n$
is irreducible and primitive, then
$\a^n/G$ is terminal whenever $n\geq 5$.
\end{say}

\section{Flips} 

For more on flips, see \cite[Chap.6]{km-book}, \cite{c-book}
or \cite{ha-mc}. 
\medskip

The following is the most general definition of flips.

\begin{defn}\label{flip.defn}  Let 
$f^-:  X^-  \to   Y$  be a proper birational morphism 
between pure dimensional $S_2$ schemes such 	that
the exceptional set $\ex (f^-)$ has codimension at least two in
$X^-$.  Let $H^-$ be an  $\r$-Cartier 
divisor on $X^-$ such that  $-H^-$
is
$f^-$-ample.
 A  pure dimensional $S_2$ scheme 
$X^+$ together with a proper birational morphism 
$f^+: X^+ \to  Y$  
	is called an {\it $H^-$-flip} of  $f^-$  if
\begin{enumerate} 
\item  the exceptional set  $\ex (f^+)$ has codimension at least
two in $X^+$.  
\item  the birational transform $H^+$ of $H^-$ on $X^+$   is $\r$-Cartier 
and $f^+$-ample.
\end{enumerate} 
By a slight   abuse of terminology, the
rational map $\phi:= \bigl(f^+)^{-1}\circ f^-: X^- \map X^+$ is also
 called an $H^-$-flip.  We will see in (\ref{ma.mu.exrc}) or
(\ref{flip=proj.exrc})
 that a flip is unique and the
main question is its existence. A flip gives the following
diagram:
$$
\begin{array}{rrcll}& X^-\!\!\! &\stackrel{\phi}{\map} &X^+&\\
\mbox{($-H^-$ is $f^-$-ample)} & f^-&\searrow \quad\swarrow &f^+ &
\mbox{($H^+$ is $f^+$-ample).}\\ &&Y&&
\end{array}
$$
\end{defn}

\begin{warning}\label{warn.flip-flop}
 In the literature the notion of flip
is  frequently used in more restrictive ways.
Here are the most commonly used variants that appear,
sometimes without explicit mention.
\begin{enumerate}
\item In older papers,  flip refers to the case
when $X^-$ is terminal  and $H=K_{X^-}$.
These are the ones needed when we start the MMP with a smooth variety.
\item 
In the MMP for  pairs  $(X,\Delta)$ we are  interested in flips when
$(X^-, \Delta^-)$ is a klt (or dlt or lc) pair 
and $H=K_{X^-}+\Delta^-$.
In older papers this is called a log-flip, but more recently
it is called simply a flip. 
\item Given $(X^-, \Delta^-)$,
a $\bigl(K_{X^-}+\Delta^-\bigr)$-flip is frequently called
a $\Delta^-$-flip.
\item The statement ``$n$-dimensional terminal (or canonical, klt,  \dots)
flips exist'' means that the $H^-$-flip of
$f^-:  X^-  \to   Y$ exists  whenever $\dim X^-=n$,
$H^-=K_{X^-}+\Delta^-$ and 
 $(X^-, \Delta^-)$ is terminal (or canonical, klt,  \dots).
\item In many cases the relative Picard number of $X^-/Y$ is 1.
Thus, up to $\r$-linear equivalence, there is a unique $f^-$-negative
divisor and the choice of $H^-$ is irrelevant; hence omitted.
This variant is frequently used for nonprojective schemes
or complex analytic spaces,
when a relatively ample divisor may not exist.
\item 
A flip is called a {\it flop} if $K_{X^-}$ is
numerically $f^-$-trivial,
or, if one has in mind a fixed 
$(X^-, \Delta^-)$, if
$K_{X^-}+\Delta^-$ is numerically $f^-$-trivial.
\item Let $X$ be a scheme and $H$ an $\r$-divisior on $X$. 
Especially when studying sequences of flips, an $H$-flip
could refer to 
 any {\it $H^-$-flip} of  $f^-:X^-\to Y$
if there is a birational  contraction
$g:X\map X^-$ and $H^-$ is the birational transform of $H$.
\end{enumerate}
\end{warning}

\begin{exrc}\label{ma.mu.exrc} Prove the following result of
Matsusaka and Mumford \cite{ma-mu}.

Let $X_i$ be pure dimensional $S_2$-schemes and
$X_i\to S$ projective morphisms with relatively ample
divisors $H_i$. Let $U_i\subset X_i$ be open subsets such that
$X_i\setminus U_i$ has codimension $\geq 2$ in $X_i$.
Let $\phi_U:U_1\to U_2$ be an isomorphism such that
$\phi_U(H_1|_{U_1})=H_2|_{U_2}$. 

Then $\phi_U$ extends to an isomorphism $\phi_X:X_1\to X_2$.
\end{exrc}

\begin{exrc} Notation as in (\ref{flip.defn}).
Prove that $f^-_*(H^-)$ is not $\r$-Cartier on $Y$.
\end{exrc}

We see in (\ref{noflip.exrc}) that not all flips exist.
Currently, the strongest existence theorem 
is the following.

\begin{thm} \cite{ha-mc, bchm} Dlt flips exist.
\end{thm}

\begin{exrc} Let $\phi:X^-\map X^+$ be a
 $(K_{X^-}+\Delta^-)$-flip.
Prove that for any $E$,
$a(E,X^-,\Delta_{X^-})\leq a(E,X^+,\Delta_{X^+})$
and the inequality is strict iff the center of $E$ on
$X^-$ is contained in $\ex(\phi)$.
\end{exrc}

\begin{defn} Let $(X,\Delta)$ be an lc pair
and $f:X\to S$ a proper morphism.
A {\it sequence of flips over $S$} starting with  $(X,\Delta)$
is a sequence of birational maps
$\phi_i$ and morphisms $f_i$ 
$$
\begin{array}{rcl} X_i\!\!\! &\stackrel{\phi_i}{\map} &X_{i+1}\\
  f_i\!\!\!&\searrow \quad\swarrow &\!\!\! f_{i+1} \\
 &S&
\end{array}
$$
(starting with
$X_0=X$) such that for every $i\geq 0$,  $\phi_i$ is a
$\bigl(K_{X_i}+\Delta_i\bigr)$-flip
where $\Delta_i$ is the birational transform of $\Delta$
on $X_i$.
\end{defn}

The basic open question in the field is the following

\begin{conj} Starting with an lc pair $(X_0,\Delta_0)$,
there is a no infinite sequence
of flips $\phi_i:(X_i,\Delta_i)\map (X_{i+1},\Delta_{i+1})$.
\end{conj}

This is known in dimension 3, almost known in dimension 4
and known in certain important cases in general;
see \cite{bchm} or (\ref{running.mmp.defn}) for more precise statements.

\begin{exrc} Let $\phi_i:(X_i,\Delta_i)\map (X_{i+1},\Delta_{i+1})$ be a 
sequence of flips. Prove that the composite
$\phi_n\circ\cdots\circ \phi_0:X_0\map X_{n+1}$
can not be an isomorphism.
\end{exrc}

\begin{prob} 
 Let $\phi_i:(X_i,\Delta_i)\map (X_{i+1},\Delta_{i+1})$ be a 
sequence of flips. Prove that $(X_n,\Delta_n) $ can not be isomorphic to 
$(X_0, \Delta_0)$
for $n>0$. (I do not know how to do this, but it may not be hard.)

By contrast, show that the involution $\tau$ in  (\ref{first.flop.exrc})
is a flop and even a flip for some
$H=K_X+\Delta$ where $(X,\Delta)$ is klt.
(Thus $X_n$ could be isomorphic to $X_0$, but the isomorphism should not
preserve $\Delta$.)
\end{prob}

\begin{exrc}[Simplest flop] Let
$L_1,L_2\subset \p^3$ be two lines intersecting at  a point $p$.
Let $X_1:=B_{L_1}B_{L_2}\p^3$ and 
$X_2:=B_{L_2}B_{L_1}\p^3$.
Set $Y:=B_{L_1+L_2}\p^3$.

Show that the identity on $\p^3$
induces morphisms $f_i:X_i\to Y$ and
a rational map
$\phi:X_1\map X_2$. We get a  flop diagram
$$
\begin{array}{rcl} X_1 &\stackrel{\phi}{\map} &X_2\\
f_1 &\searrow \quad\swarrow & f_2\\ 
&Y.&
\end{array}
$$
 Show that neither $\phi$ nor $\phi^{-1}$
contracts  divisors but neither is a morphism.
Describe how to factor $\phi$ into a composite of
smooth blow ups and blow downs.
\end{exrc}

\begin{exrc}[Non-algebraic flops] Let $X\subset \p^4$ be a general  smooth
quintic hypersurface. It is know that for every $d\geq 1$,
$X$ contains a smooth rational curve $\p^1\cong C_d\subset X$
with normal bundle $\o_{\p^1}(-1)+\o_{\p^1}(-1)$
\cite{clem}. 

Prove that the flop of $C_d$ exists if we work with
compact complex manifolds. Denote the flop by
 $\phi_d:X\map X_d$  and let
$H_d\in H^2(X_d,\z)$ be the image of the hyperplane class.
Compute the self-intersection $(H_d^3)$.
Conclude that the $X_d$ are not homeomorphic to each other and
not projective.
\end{exrc}

\begin{exrc}[Harder flops] Let
$C_1,C_2\subset \p^3$ be two smooth curves intersecting at  a single 
point $p$ where they are tangent to order $m$.
Let $X_1:=B_{C_1}B_{C_2}\p^3$ and 
$X_2:=B_{C_2}B_{C_1}\p^3$.
Set $Y:=B_{C_1+C_2}\p^3$.

Show that the identity on $\p^3$
induces morphisms $f_i:X_i\to Y$,
a rational map
$\phi:X_1\map X_2$ and we get a  flop diagram
as before.
Describe how to factor $f$ into a composite of
smooth blow ups and blow downs.
\end{exrc}

\begin{exrc}[Even harder flops]
Consider the variety
$$
X:=(sx+ty+uz=sz^2+tx^2+uy^2=0)\subset \p^2_{xyz}\times \a^3_{stu}.
$$
Show that $X$ is smooth, the projection
$\pi:X\to \a^3$ has degree 2 and $C:=\red\pi^{-1}(0,0,0)$
is a smooth rational curve.
Compute $(C\cdot K_X)$ and the normal bundle of $C$.

Let $Y\to \a^3$ be the normalization of $\a^3$ in $k(X)$.
Determine the singularity of $Y$ sitting over the origin.

As before, the Galois involution of $Y\to \a^3$
provides the flop of $X\to Y$.

It is quite tricky to factor $f$ into a composite of
smooth blow ups and blow downs.
\end{exrc}

\begin{exrc}[Simplest flips]\label{simp.flip.exrc}
 Fix $n\geq 3$ and consider the affine hypersurface
$$
Z:=(u^n-u^{n-1}y+x^{n-1}z=0)\subset \a^4,
$$
which we view as a degree $n$ covering of the
$(x,y,z)$-space.

Show that $Z$ is not normal and its normalization has a
unique singular point which lies above $(0,0,0)$.

Show that
$$
X^+:=(s^nx-s^{n-1}ty+t^nx^{n-1}z=0)\subset \a^3_{xyz}\times \p^1_{st}
$$
is a small resolution of  $Z$. Write down the morphism
$X^+\to Z$. It has a unique 1-dimensional fiber
$C^+\subset X^+$. Determine the normal bundle
of $C^+$ in $X^+$ and the intersection number of
$C^+$ with the canonical class.

Construct another small modification
$X^-\to Z$ as follows.
First blow up the ideal $(z, u^{n-1})$.
We get the variety $X_1$ defined by equations
$$
\bigl(s(y-u)-tx^{n-1}= sz-tu^{n-1}= u^n-u^{n-1}y+x^{n-1}z=0\bigr)
\subset \a^4_{xyzu}\times \p^1_{st}.
$$
Show that  the $s\neq 0$ chart is smooth
and on the $t\neq 0$ chart we have a complete intersection
$$
\bigl(w(y-u)-x^{n-1}= wz-u^{n-1}=0\bigr)
\subset \a^4_{xyzuw}\qtq{with $w=s/t$.}
$$
Setting $y':=y-u$ we have the local  equations for $X_1$
$$
wy'-x^{n-1}= wz-u^{n-1}=0.
$$
Write down a $\z/(n-1)$-invariant finite morphism
to the above local chart on $X_1$
 from $\a^3_{pqr}$ with the $\z/(n-1)$-action
$(p,q,r)\mapsto (\epsilon p,\epsilon q,\epsilon^{-1}r)$, 
where $\epsilon$ is a primitive $(n-1)$-st root of unity.
Let $X^-$ be the normalization of $X_1$.
Show that $X^-$ has a single quotient singularity of
the above form. 

 Write down the morphism
$X^-\to Z$. It has a unique 1-dimensional fiber
$C^-\subset X^-$. Determine  the intersection number of
$C^-$ with the canonical class.
\end{exrc}

\begin{exrc}
Let now $Y$ be any smooth 3-fold and $L$ a very ample
line bundle on $Y$ with 3 general sections $f,g,h$.
 Fix $n\geq 3$ and consider the  hypersurface
$$
Z:=(u^n-u^{n-1}g+f^{n-1}h=0)\subset L^{-1}.
$$
One small resolution is given by
$$
X^+:=(s^nf-s^{n-1}tg+t^nh=0)\subset Y\times \p^1_{st}.
$$
Compute its canonical class in terms of
 $K_Y$ and $L$. 
\end{exrc}

\begin{exrc}[Log terminal flips]
 Work out the analog of (\ref{simp.flip.exrc}) 
when we start with
$$
X^+:=(s^nx-s^{n-i}t^iy+t^nz=0)\subset \a^3_{xyz}\times \p^1_{st}.
$$
\end{exrc}

\begin{exrc}\label{flip=proj.exrc}
 Let $X$ be a Noetherian, reduced, pure dimensional, $S_2$-scheme
and $D$ a Weil divisor on $X$
which is Cartier in codimension 1. Prove that the following are equivalent.
\begin{enumerate}
\item $\sum_{m\geq 0} \o_X(mD)$ is a finitely generated sheaf of 
$\o_X$-algebras.
\item There is a proper, birational morphism
$\pi:X^+\to X$ such that the exceptional set $\ex(\pi)$ has codimension
$\geq 2$ and the birational transform $D^+:=\pi^{-1}_*(D)$ is
$\q$-Cartier and $\pi$-ample.
\end{enumerate}

Hint of proof. (2) $\Rightarrow$ (1) is easy.

To see the converse, set $X^+:=\proj_X\sum_{m\geq 0} \o_X(mD)$.
We need to show that  $X^+\to X$ is small.
Assume that $E\subset \ex(\pi)$ is an exceptional divisor.
Study the sequence
$$
0\to \o_{X^+}(mD^+)\to \o_{X^+}(mD^++E)\to\o_{E}\bigl((mD^++E)|_E\bigr)\to 0
$$
to get, for some $m>0$,  a section of $\o_{X^+}(mD^++E)$
which is not a section of  $\o_{X^+}(mD^+)$.
By pushing forward to $X$, we would get extra sections
of $\o_X(mD)$.
\end{exrc}

\begin{exrc}\label{flip=fg.exrc}
Let $(X,\Delta)$ be klt. Let
$f:X\to Y$ be a small $(K_X+\Delta)$-negative contraction.
Show that there is a $\q$-divisor $D$ on $X$ such that
 $(X,\Delta+D)$ is klt
and $(K_X+\Delta+D)\sim_{\q,f} 0$.

Conclude from this that  $\bigl(Y, f_*(\Delta+D)\bigr)$
 is klt.
\end{exrc}

A consequence of the relative MMP is the
following finite generation result, which we prove in (\ref{D.amp.qf}).
By (\ref{flip=fg.exrc}), it formally implies the existence of dlt flips.

\begin{thm}\label{klt.fg.thm}
 Let $(X,\Delta)$ be klt and  $D$  a $\q$-divisor on $X$.
 Then  $\sum_{m\geq 0} \o_X(\rdown{mD})$ is a finitely generated sheaf of 
$\o_X$-algebras.
\end{thm}

\begin{exrc}  Show that $\rdown{A+B}\geq \rdown{A}+\rdown{B}$
for any divisors $A,B$, thus,
for any divisor $D$, 
$R(X,D):=\sum_{m\geq 0} H^0(X,\o_X(\rdown{mD}))$ is a ring.

Give examples where
$R^u(X,D):=\sum_{m\geq 0} H^0(X, \o_X(\rup{mD}))$
 is not a ring. Note, however, that  
$\rup{A+B}\geq \rdown{A}+\rup{B}$, thus
$R^u(X,D)$ is an $R(X,D)$-module.
\end{exrc}

\begin{exrc}  Let $X$ be normal and  $D$ an  $\r$-divisor.
Show that if 
 $\sum_{m\geq 0} \o_X(\rdown{mD})$ is a finitely generated sheaf of 
$\o_X$-algebras then $D$ is  a $\q$-divisor.
\end{exrc}

The following example shows that (\ref{klt.fg.thm})
fails for lc pairs. 

\begin{exrc}\label{nonfg.exrc}
 Let $E\subset \p^2$ be a smooth cubic.
Let $S$ be obtained by blowing up $9$ general points
on $E$ and let $E_S\subset S$ be the birational transform of $E$.
Let $H$ be a sufficiently ample divisor on $S$
giving a projectively normal embedding $S\subset \p^n$.
Let $X\subset \a^{n+1}$ be the cone over $S$ and $D\subset X$
the cone over $E_S$.

Prove that $(X,D)$ is lc yet
 $\sum_{m\geq 0} \o_X(mD)$ is not a finitely generated sheaf of 
$\o_X$-algebras.
\medskip

Hints.  First show that
$H^0(X, \o_X(mD))=\sum_{r\geq 0} H^0\bigl(S, \o_S(mE_S+rH)\bigr)$.
Check that $\o_S(mE_S+rH)$ is very ample if $r>0$
but $\o_S(mE_S)$ has only the obvious section which vanishes along $mE_S$.
Thus the multiplication maps
$$
\sum_{a= 0}^{m-1} 
H^0\bigl(S, \o_S(aE_S+H)\bigr)\otimes H^0\bigl(S, \o_S((m-a)E_S)\bigr)
\to H^0\bigl(S, \o_S(mE_S+H)\bigr)
$$
are never surjective.
\end{exrc}

The next exercise shows that log canonical flops
sometimes do not exist.

\begin{exrc}  \label{noflip.exrc}
Let $E$ be an elliptic curve, $L$ a degree 0 non-torsion line bundle
and $S=\p_E(\o_E+L)$. Let $C_1, C_2\subset S$ be the corresponding
sections of $S\to E$.
Note that $K_S+C_1+C_2\sim 0$.
Let $0\in X$ be a cone over $S$ and $D_i\subset X$ the cones over $C_i$.
Show that $(X,D_1+D_2)$ is lc. 

Following the method of
(\ref{nonfg.exrc}) show that 
$\sum_{m\geq 0} \o_X(mD_i)$ is not a finitely generated sheaf of 
$\o_X$-algebras for $i=1,2$.

Let $F\subset S$ be a fiber of $S\to E$ and
$B\subset X$ the cone over $F$. Show that
$\sum_{m\geq 0} \o_X(mB)$ is  a finitely generated sheaf of 
$\o_X$-algebras and describe the corresponding
small contraction $\pi:Z\to X$. 

Prove that the flip of $\pi:Z\to X$ does not exist
(no matter what $H$ we choose). 

What happens if $L$ is  a torsion element in $\pic(E)$?
\end{exrc}

\begin{exrc} Let $S$ be a Noetherian, reduced, 2-dimensional, $S_2$-scheme
and $D$ a Weil divisor on $S$. Prove that 
 $\sum_{m\geq 0} \o_S(mD)$ is a finitely generated sheaf of 
$\o_S$-algebras iff $\o_S(mD)$ is locally free for some $m>0$.

Use this to show that the following algebras are not finitely generated.
\begin{enumerate}
\item $S$ is a cone over an elliptic curve and $D\subset S$ a general line.
State the precise generality condition.
\item Let $C\subset \p^n$ be a projectively normal curve of genus $\geq 2$
and $S\subset \a^{n+1}$ the cone over $C$.  Assume that
$\o_C(1)$ is a general line bundle and let $D=K_S$.
Again, state the precise generality condition.
\item Let $S$ be the quadric cone $(xy-z^2=0)\subset \a^3$ and the $(u,v)$-plane
glued together along the lines $(x=z=0)$ and $(v=0)$.
(Show that this surface does not embed in $\a^3$ but 
realize it in $\a^4$ by explicit equations.)
Set $D=K_S$.
\end{enumerate}
\end{exrc}

The following conjecture is known if $x\in H$ is a quotient singularity
\cite{ksb} or when $x\in H$ is a quadruple point \cite{stev}.
It is quite remarkable that, aside from the case
when $x\in H$ is a quotient singularity, the conjecture
seems unrelated to the minimal model program.

\begin{conj}\cite[6.2.1]{k-ffmm} Let $x\in X$ be a 3-dimensional
normal singularity and $x\in H\subset X$ a Cartier divisor.
Assume that $x\in H$ is a (normal) rational surface singularity.
Then $\sum_{m\geq 0} \o_X(mK_X)$ is a finitely generated sheaf of 
$\o_X$-algebras. 
\end{conj}

\section{Minimal models}

For more details, see \cite[3.7--8]{km-book}
or \cite{bchm}.
\medskip

\begin{defn}[Running the MMP] \label{running.mmp.defn}
Let $(X,\Delta)$ be a   
pair such that $K_X+\Delta$ is $\q$-Cartier
and $f:X\to S$ a proper morphism. 
Assume for simplicity that $X$ is $\q$-factorial.
A {\it running of the $(K_X+\Delta)$-MMP} over $S$ yields a
sequence
$$
(X,\Delta)=:(X_0,\Delta_0)\stackrel{\phi_0}{\map}
(X_1,\Delta_1)\stackrel{\phi_1}{\map}\quad\cdots\quad
\stackrel{\phi_{n-1}}{\map}
(X_r,\Delta_r),
$$
where each $\phi_i$ is either the divisorial
contraction of a $(K_{X_i}+\Delta_i)$-negative extremal ray
or the flip of a small 
contraction of a $(K_{X_i}+\Delta_i)$-negative extremal ray, 
$\Delta_{i+1}:=(\phi_i)_*\Delta_i$ and all the $X_i$ are
$S$-schemes $f_i:X_i\to S$ such that $f_i=f_{i+1}\circ \phi_i$.
We say the the $(K_X+\Delta)$-MMP {\it stops} or {\it terminates}
with $(X_r,\Delta_r)$ if 
\begin{enumerate}
\item either $K_{X_r}+\Delta_r$ is $f_r$-nef
(and there are no more extremal rays),
\item or there is a Fano contraction $X_r\to Z_r$.
\end{enumerate}
Sometimes we impose a stronger restriction:
\begin{enumerate}\setcounter{enumi}{2}
\item[(2')] every extremal contraction of  $(X_r,\Delta_r)$
is Fano.
\end{enumerate}

Conjecturally, every running of the $(K_X+\Delta)$-MMP stops.
This is known if $\dim X\leq 3$  \ \cite{kaw-term},
in many cases in dimension 4 \cite{ahk}
or when the generic fiber of $f$ is of general type \cite{bchm}
and at each step the extremal rays are chosen ``suitably.''
Note that the latter includes the case when $f$ is birational
(or generically finite),
since a point is a 0-dimensional variety of general type.

(Everything works the same if $X$ is not $\q$-factorial, except
in that case it does not make sense to distinguish
divisorial contractions and flips.) 
\end{defn}

\begin{defn} \label{min.defn} Let $(X,\Delta)$ be a  pair and $f:X\to S$ a
proper morphism.  We say that $(X,\Delta)$ is an
$$
\left.
\begin{array}{l}
\mbox{$f$-weak canonical} \\
\mbox{$f$-canonical} \\
\mbox{$f$-minimal} 
\end{array}
\right\}
\mbox{ model if $(X,\Delta)$ is }
\left\{
\begin{array}{c}
\mbox{lc} \\
\mbox{lc} \\
\mbox{dlt} 
\end{array}
\right\}
\mbox{ and $K_X+\Delta$ is }
\left\{
\begin{array}{c}
\mbox{$f$-nef} \\
\mbox{$f$-ample} \\
\mbox{$f$-nef} 
\end{array}
\right\}.
$$
\end{defn}

\begin{warning} Note that a canonical model $(X,\Delta)$
has {\em log} canonical singularities, not necessarily
canonical singularities. This, by now entrenched, unfortunate
terminology is a result of an incomplete shift.
Originally everything was defined only  for $\Delta=0$.
When $\Delta$ was introduced, its presence was
indicated by putting ``log'' in front of adjectives.
Later, when the use of $\Delta$ became pervasive,
people started dropping the prefix ``log''.
This is usually not a problem. For instance,
the canonical ring $R(X,K_X)$ is just the
$\Delta=0$ special case of the log canonical ring $R(X,K_X+\Delta)$.

However, canonical singularities are not the 
$\Delta=0$ special cases of log canonical singularities.
\end{warning}

\begin{defn}\label{3-6.minmod}  Let $(X,\Delta)$ be a   
pair such that $K_X+\Delta$ is $\q$-Cartier
and $f:X\to S$ a proper morphism.   A pair
$(X^w,\Delta^w)$ sitting in a diagram
$$
\begin{array}{rcl} X &\stackrel{\phi}{\map} &X^w\\ f
&\searrow \quad\swarrow & f^w\\  &S&
\end{array}
$$
 is called a {\it  weak canonical model  of}
$(X,\Delta)$ over $S$ if 
\begin{enumerate} 

\item $f^w$ is proper, 

\item $\phi$ is a contraction, that is, 
$\phi^{-1}$ has no exceptional divisors, 

\item  $\Delta^w=\phi_*\Delta$, 

\item $K_{X^w}+\Delta^w$ is $\q$-Cartier and $f^w$-nef,   and  

\item  $a(E,X,\Delta)\leq a(E,X^w,\Delta^w)$ for every
$\phi$-exceptional divisor $E\subset X$. 
Equivalently,
$(K_{X}+\Delta)-\phi^*\bigl(K_{X^w}+\Delta^w\bigr)$
 is effective and $\phi$-exceptional.
\end{enumerate}

A  weak canonical model      
$(X^m,\Delta^m)=(X^w,\Delta^w)$ is called a {\it  minimal
model of}
 $(X,\Delta)$ over $S$ if,  in addition to (1--4), we have 
\begin{enumerate} 
\item[($5^m$)]    $a(E,X,\Delta)< a(E,X^m,\Delta^m)$ for
every $\phi$-exceptional divisor $E\subset X$.
\end{enumerate}

A weak canonical model 
$(X^c,\Delta^c)=(X^w,\Delta^w)$ is called a {\it 
canonical model  of}
$(X,\Delta)$ over $S$ if,  in addition to (1--3) and (5) we have 
\begin{enumerate} 
\item[($4^c$)] $K_{X^c}+\Delta^c$ is  $\q$-Cartier and $f^c$-ample.
\end{enumerate}
\end{defn}

\begin{exrc} 
Let $(X,\Delta)$ be a   
pair such that $K_X+\Delta$ is $\q$-Cartier
and $f:X\to S$ a proper morphism. Run the MMP:
$$
(X,\Delta)=:(X_0,\Delta_0)\stackrel{\phi_0}{\map}
(X_1,\Delta_1)\stackrel{\phi_1}{\map}\quad\cdots\quad
\stackrel{\phi_{n-1}}{\map}
(X_r,\Delta_r),
$$
and assume that $K_{X_r}+\Delta_r$ is $f$-nef.
Show that $(X_r,\Delta_r)$ is a minimal model of
$(X,\Delta)$ over $S$.
\end{exrc}

\begin{exrc} Let $f:(X,\Delta)\to S$ be a canonical model.
Let $g:X'\to X$ be a proper birational morphism
with exceptional divisors $E_i$. When is
$f:(X,\Delta)\to S$  a canonical model of
$(X', g^{-1}_*\Delta+\sum e_iE_i)$?
\end{exrc}

\begin{exrc} \label{3-6.minmod.sings}  Let
$\phi:(X,\Delta)\map (X^w,\Delta^w)$ be a   weak canonical model. 
Prove that 
$$
a(E,X^w,\Delta^w)\geq a(E,X,\Delta) 
\qtq{ for every divisor $E$.} 
$$

Hint.  Fix $E$ and consider any diagram
$$
\begin{array}{rcl} &Y&\\ g &\swarrow \quad\searrow & h\\ 
 X &\stackrel{\phi}{\map} &X^w\\ f &\searrow \quad\swarrow
& f^w\\  &S&
\end{array}
$$ where  $\cent_YE$ is a  divisor.  Write
$K_Y$ in two different ways and apply (\ref{exc.eff.exrc}).
\end{exrc} 

\begin{exrc} \label{3-6.minmod.int.no}  Let
$\phi:(X,\Delta)\map (X^w,\Delta^w)$ be a   weak canonical model. 
Prove that if a curve  $C\subset X$ is not contained
in $\ex(\phi)$ then
$$
C\cdot (K_X+\Delta)\geq \phi_*(C)\cdot (K_{X^w}+\Delta^w).
$$
\end{exrc}

\begin{exrc}\label{exc.eff.exrc} Let $h:Z\to Y$ be a proper
birational morphism between normal varieties. Let
$-B$ be an $h$-nef $\q$-Cartier $\q$-divisor on
$Z$. Then
\begin{enumerate}

\item $B$ is effective iff $h_*B$ is.

\item Assume that $B$ is effective. Then for every $y\in Y$,
either $h^{-1}(y)\subset\supp B$ or
$h^{-1}(y)\cap\supp B=\emptyset$.
\end{enumerate}

Hint.  Use induction on $\dim Z$ by passing to a
hyperplane section  $H\subset Z$. Be careful:
 $h_*(B\cap H)$ need not be contained in  $h_*B$.
\end{exrc}

\begin{exrc}[$\q$-factorialization] \label{q-fact.exrc}
Let $(X,\Delta)$ be klt. Let $f:Y\to X$ be a log resolution
with exceptional divisor $E$. For $0<\epsilon\ll 1$
run the $\bigl(Y, f^{-1}_*\Delta+(1-\epsilon)E\bigr)$-MMP over $X$ and
assume that it stops.
(This is not a restriction by (\ref{running.mmp.defn}).)

Prove that the MMP stops at a small contraction $f_r:Y_r\to X$
such that $Y_r$ is $\q$-factorial.

It is called a {\it $\q$-factorialization} of $X$.

More generally, prove that  $\q$-factorializations
exist if  $(X,\Delta)$ is dlt. Find lc examples
without any $\q$-factorialization.
\end{exrc}

\begin{exrc} \label{D.amp.qf}
Notation as in (\ref{q-fact.exrc}). Let $D$ be any
Weil divisor on $X$. Prove that there is a
 $\q$-factorialization $f_D:Y_D\to X$ such that
the birational transform of $D$ on $Y_D$ is $f_D$-nef.

Use this to prove that  $\q$-factorializations are never unique,
save when $X$ itself is $\q$-factorial.

Use this and the contraction theorem to prove (\ref{klt.fg.thm}).
\end{exrc}

\begin{warning} You may have noticed already that we have not defined
when a pair $(X',\Delta')$ is birational to another pair $(X,\Delta)$.
The problem is: what should  the coefficient of a
divisor $D\subset X'$  be in $\Delta'$ when the center of
$D$ on $X$ is not a divisor. 

One approach is to insist that birational pairs
have the same canonical rings. Then
the next exercise suggests a definition. 

It is, however, best to keep in mind that
birational equivalence of pairs is a problematic concept.
\end{warning}

\begin{exrc} \label{can.ring.=.cond}
Let $f_1:X_1\to S$ and $f_2:X_2\to S$ be proper morphisms
of normal schemes
and  $\phi:X_1\map X_2$ a birational map such that $f_1=f_2\circ \phi$.
Let $\Delta_1$ and $\Delta_2$ be $\q$-divisors 
such that $K_{X_1}+\Delta_1$ and  $K_{X_2}+\Delta_2$ are $\q$-Cartier.
Prove that 
$$
{f_1}_*\o_{X_1}(mK_{X_1}+\rdown{m\Delta_1})=
 {f_2}_*\o_{X_2}(mK_{X_2}+\rdown{m\Delta_2})
\qtq{for $m\geq 0$}
$$
if the following conditions hold:
\begin{enumerate}
\item $a(E,X_1,\Delta_1)=a(E,X_2,\Delta_2)$ if 
$\phi$ is a local isomorphism at the generic point of $E$,
\item $a(E,X_1,\Delta_1)\leq a(E,X_2,\Delta_2)$ if 
$E\subset X_1$ is $\phi$-exceptional, and
\item $a(E,X_1,\Delta_1)\geq a(E,X_2,\Delta_2)$ if 
$E\subset X_2$ is $\phi^{-1}$-exceptional.
\end{enumerate}

Hints:  Let $Y$ be the normalization of the closed graph of $\phi$
in $X_1\times_SX_2$ and
$g_i:Y\to X_i$ the projections. We can write
$$
\begin{array}{rcl}
K_Y&\simq& g_1^*(K_{X_1}+\Delta_1)+\sum_E a(E,X_1,\Delta_1)E,\qtq{and}\\
K_Y&\simq& g_1^*(K_{X_2}+\Delta_2)+\sum_E a(E,X_2,\Delta_2)E.
\end{array}
$$
Set $b(E):=\max\{-a(E,X_1,\Delta_1), -a(E,X_2,\Delta_2)\}$.
Prove that  $\sum_E \bigl(b(E)+a(E,X_i,\Delta_i)\bigr)E$ 
is effective and $g_i$-exceptional for $i=1,2$.
Conclude that
$$
\begin{array}{l}
{(f_i\circ g_i)}_*\o_{Y}(mK_Y+\textstyle{\sum_E}mb(E)E)\\
\qquad = {f_i}_* {g_i}_*
\o_{Y}\Bigl(g_1^*(mK_{X_1}+m\Delta_1)+
\sum_E \bigl(mb(E)+ma(E,X_1,\Delta_1)\bigr)E\Bigr)\\
\qquad ={f_i}_* \o_{X_i}(mK_{X_i}+m\Delta_i).
\end{array}
$$
\end{exrc}

\begin{exrc}\label{3-6.canmod.unique}  Let $(X,\Delta)$ be
a  lc pair with $\Delta\geq 0$,  
$f:X\to S$ a proper morphism
and $f^w: (X^w,\Delta^w)\to S$  a weak minimal model. Prove the
following:  
\begin{enumerate}
\item $f_*\o_X(mK_X+\rdown{m\Delta})=f^w_*\o_{X^w}(mK_{X^w}+\rdown{m\Delta^w})$
for every $m\geq 0$.
\item If a canonical model    $(X^c,\Delta^c)$ exists then
$$
X^c=\proj_S\textstyle{\sum}_{m\geq 0} f_*\o_X(mK_X+\rdown{m\Delta}),
$$
and the right hand side is a sheaf of  finitely generated algebras.
In particular,  a canonical model  is unique.
\item Any two minimal models of $(X,\Delta)$ are isomorphic
in codimension one. (Hint: Prove this first when
$\Delta=0$ and $(X,0)$ is terminal. The general case is
more subtle.)
\end{enumerate}
\end{exrc}

\begin{exrc} Assume  that $X$ is irreducible,
$$
R(X,K_X+\Delta):= \textstyle{\sum}_{m\geq 0} f_*\o_X(mK_X+\rdown{m\Delta})
$$
 is a sheaf of  finitely generated algebras and
$$
\dim X=\dim \proj_SR(X,K_X+\Delta).
$$
Prove that the natural map
$\phi:X\map  \proj_SR(X,K_X+\Delta)$ is birational
and 
$$
(X^c,\Delta^c):=\bigl(\proj_SR(X,K_X+\Delta), \phi_*\Delta\bigr)
$$
is the canonical model of $(X,\Delta)$.
\medskip

Hint: You should find  (\ref{fg.bir.maps.lem}) useful.
\end{exrc}

\begin{exrc}\label{fg.bir.maps.lem}
 Let $X$ be an irreducible and normal scheme, $L$ a Weil divisor
on $X$ and $f:X\to S$ a proper morphism, $S$ affine.
Write $|L|=|M|+F$ where $|M|$ is the moving part and
$F$ the fixed part.
Assume that $R(X,L):=\textstyle{\sum}_{m\geq 0} f_*\o_X(mL)$
is  generated by $f_*\o_X(L)$.
Set $Z:=\proj_S R(X,L)$ with projection $p:Z\to S$ and let 
$\phi: X\map Z$
be the natural morphism. Prove that
\begin{enumerate}
\item $Z\setminus \phi(X)$
has codimension $\geq 2$ in $Z$.
\item If $\phi$ is generically finite then it is  birational  and
  $F$ is $\phi$-exceptional.
\end{enumerate}
(Hint: This is similar to (\ref{flip=proj.exrc}).)
\end{exrc}

\begin{exrc}[Chambers in the cone of  big divisors] 
Let $X$ be a normal variety and
$D_i$ big $\q$-divisors. Assume that the rings
$$
R(D_i):=\tsum_{m\geq 0}H^0\bigl(X, \o_X(\rdown{mD_i})\bigr)
$$
are finitely generated and the maps
$X\map \proj R(D_i)$ are birational and independent of $i$. Let
$D=\sum a_iD_i$ be a nonnegative $\q$-linear combination.

Prove that   $R(D)$
is finitely generated and 
$X\map \proj R(D)$ is the same map as before.  

Conclude that the set of all big $\q$-divisors with the same
$X\map \proj R(D)$ forms a convex subcone, called a
{\it chamber} in the cone of  big divisors.
\end{exrc}

\begin{exrc} Develop a relative version of the
notion of chambers of divisors for maps.
(Note that for birational maps, every divisor is relatively big.)

Let $Y\to X$ be  a $\q$-factorialization of a klt pair
$(X,\Delta)$ (\ref{q-fact.exrc}).
 Prove that there is a one-to-one correspondence
between open chambers of $N^1(Y/X)$ and
$\q$-factorializations of $X$.

What kind of maps correspond to the other chambers?
\end{exrc}

\begin{exrc} Let $a_i$ be different complex numbers.
Consider the singularity
$$
X=X(a_1,\dots,a_n):=\bigl(xy-\prod_i(u-a_iv)=0\bigr)\subset \a^4.
$$
Find a small resolution of $X$ by repeatedly blowing up
planes of the form $(x=u-a_iv=0)$.

Prove that the class group  $\cl(X)$ of $X$ is  generated by the planes
$(x=u-a_iv=0)$, with a single relation
$\sum_i [x=u-a_iv=0]=0$.

Describe all small resolutions of $X$ and the
corresponding chamber structure on $\cl(X)$.

(The same method can be used to describe the class group
and the chamber structure for any $cA$-type terminal
3-fold singularity, see \cite[2.2.7]{k-ffmm}. A similarly explicit
description is not known for the $cD$ and $cE$-type cases.)
\end{exrc}

\begin{exrc} Let $S:=(xy-z^3=0)\subset \a^3$
and $f:X\to S$ its minimal resolution
with exceptional curves $D_1, D_2$. Let $D_3, D_4$
be the birational transforms of the lines $(x=z=0)$
and $(y=z=0)$. For $0\leq a_i\leq 1$
describe  minimal and canonical models
of $(X,\sum a_iD_i)$ over $S$. Describe the
chamber decomposition of $[0,1]^4$.
\end{exrc}

\begin{exrc} Let $S$ be one of the singularities in
(\ref{ell.lowmult.quot.list}) and $f:X\to S$ its minimal resolution
with exceptional curves $D_i$. 
 For $0\leq a_i\leq 1$
describe  minimal and canonical models
of $(X,\sum a_iD_i)$ over $S$ and the corresponding
chamber decomposition.

(This is pretty easy for the $\z/2$-quotient.
Some of the others have many curves to check.)
\end{exrc}

For the  theory behind the next exercises, see \cite{ko-la}.

\begin{exrc} 
Let $E$ be the projective elliptic curve with affine
equation $(y^2=x^3-1)$ and set 
$\tau:(x,y)\mapsto (x,-y)$. Check that
\begin{enumerate}
\item $E/\tau\cong \p^1$.
\item $(E\times E)/(\tau\times \tau)$ has Kodaira dimension 0.
It is an example of a Kummer surface.
If $u=y_1y_2$ then it has affine equation
$$
u^2=(x_1^3-1)(x_2^3-1).
$$
Find the singularities using this equation.
\item For $n\geq 3$, $\bigl(E^n\bigr)/(\tau,\dots,\tau)$
has Kodaira dimension 0.
\end{enumerate}
\end{exrc}

\begin{exrc} 
Let $E$ be the projective elliptic curve with affine
equation $(y^3=x^3-1)$ and set 
$\sigma:(x,y)\mapsto (x,\epsilon y)$ where $\epsilon=\sqrt[3]{1}$. Check that
\begin{enumerate}
\item $E/\sigma\cong \p^1$.
\item $(E\times E)/(\sigma, \sigma^2)$ has Kodaira dimension 0.
It is an example of a K3 surface.
If $u=y_1y_2$ then it has affine equation
$$
u^3=(x_1^3-1)(x_2^3-1).
$$
Find the singularities using this equation.
\item $(E\times E)/(\sigma, \sigma)$.
If $v=y_1y_2^2$ then it has affine equation
$$
v^3=(x_1^3-1)(x_2^3-1)^2.
$$
Find the singularities using this equation.

Prove that this surface is rational in two ways:
\begin{enumerate}
\item Find many rational curves on it as preimages
of rational curves of bi-degree $(2,2)$ on $\p^1\times \p^1$.
\item Show that it is birational (even over $\z$)
to the cubic surface $y_1^3-y_2^3=x_1^3-1$.
\end{enumerate}
\item For $n\geq 3$, $\bigl(E^n\bigr)/(\sigma,\dots,\sigma)$
has Kodaira dimension 0.
\end{enumerate}
\end{exrc}

\begin{exrc} 
Let $E$ be the projective elliptic curve with affine
equation $\bigl(y^6=x(x-1)^2(x+1)^3\bigr)$ and set 
$\rho:(x,y)\mapsto (x,\epsilon y)$ where $\epsilon=\sqrt[6]{1}$. Check that
\begin{enumerate}
\item $E/\rho\cong \p^1$.
\item For $2\leq n\leq 5$, $\bigl(E^n\bigr)/(\rho,\dots,\rho)$
is uniruled, that is, it has a covering family of  rational curves.
   Try to find explicitly such a family.
(Such a family exists by \cite{ko-la}, but I do not know how to
construct one.)
I don't know if these examples are rational or unirational.
 \item For $6\leq n$, $\bigl(E^n\bigr)/(\rho,\dots,\rho)$
has Kodaira dimension 0. 
\end{enumerate}
\end{exrc}

 \begin{ack}  I thank Ch.\ Hacon, S.\ Mori and Ch.\ Xu
and the  participants of the OSU workshop organized by H.\ Clemens 
for many useful comments and corrections. 
Partial financial support  was provided by  the NSF under grant number 
DMS-0500198. 
\end{ack}

\bibliography{refs}

\vskip1cm

\noindent Princeton University, Princeton NJ 08544-1000

\begin{verbatim}kollar@math.princeton.edu\end{verbatim}

\end{document}